\documentclass[review,hidelinks,onefignum,onetabnum]{siamart220329}
\pdfoutput=1

\usepackage{rotating}
\usepackage{subfig}
\usepackage{amsmath}
\usepackage[english]{babel}
\usepackage{setspace}
\usepackage{amsfonts}
\usepackage{graphicx}
\usepackage{epstopdf}
\usepackage{algorithmic}
\ifpdf
  \DeclareGraphicsExtensions{.eps,.pdf,.png,.jpg}
\else
  \DeclareGraphicsExtensions{.eps}
\fi

\newsiamremark{remark}{Remark}
\newsiamremark{hypothesis}{Hypothesis}
\crefname{hypothesis}{Hypothesis}{Hypotheses}
\newsiamthm{claim}{Claim}

\headers{SPECTRAL APPROXIMATION OF FREDHOLM CONVOLUTION OPERATORS}{}

\title{SPECTRAL APPROXIMATION OF CONVOLUTION OPERATORS OF FREDHOLM TYPE\thanks{Submitted to the editors DATE.}}

\author{Xiaolin Liu\thanks{School of Mathematical Sciences, University of Science and Technology of China, 96 Jinzhai Road, Hefei 230026, Anhui, China (\email{xiaolin@mail.ustc.edu.cn}, \email{dengkuan@mail.ustc.edu.cn}, \email{kuanxu@ustc.edu.cn}).} \and Kuan Deng\footnotemark[2] \and Kuan Xu\footnotemark[2]}
\usepackage{amsopn}
\DeclareMathOperator{\diag}{diag}

\def\mF{\mathcal{F}}

\def\mO{\mathcal{O}}

\def\mV{\mathcal{V}}

\def\md{\mathrm{d}}

\def\bf{\boldsymbol{f}}

\def\diag{\mathop{\mathrm{diag}}}

\graphicspath{{figures/}}

\begin{document}
%

\maketitle
\begin{abstract}
We have developed a method for constructing spectral approximations for convolution operators of Fredholm type. The algorithm we propose is numerically stable and takes advantage of the recurrence relations satisfied by the entries of such a matrix approximation. When used for computing the Fredholm convolution of two given functions, such approximations produce the convolution more rapidly than the state-of-the-art methods. The proposed approximation also leads to a spectral method for solving the Fredholm convolution integral equations and enables the computation of eigenvalues and pseudospectra of Fredholm convolution operators, which is otherwise intractable with existing techniques.
\end{abstract}



\section{Introduction}
Convolution is a fundamental operation which abounds in science and engineering. The convolution of two smooth signals or functions typically results in a piecewise smooth function, which consists of three pieces. The two end pieces are known as Volterra convolution integral because each has a variable integration limit, while the middle piece, which has constant limits, is a Fredholm convolution integral \cite{hal2,xu1}. In this paper, we are concerned with the construction of approximations to the convolution operators of Fredholm type and the applications of such approximations, including the calculation of Fredholm convolution integrals and the spectral method for solving Fredholm convolution integral equations. This paper can be considered a companion work to \cite{xu2}, which addresses the spectral approximation of the convolution operators of Volterra type.

Let $f(x):[a,b] \rightarrow \mathbb{R}$ and $g(x):[c,d] \rightarrow \mathbb{R}$ be two continuous integrable compactly supported functions with $b-a>d-c$. The Fredholm convolution of $f(x)$ and $g(x)$ is
\begin{align*}
h(x) = \int_{c}^{d} f(x-t)g(t) \md t, ~~~  x \in [a+d, b+c].
\end{align*}
Since it is convenient to discuss in terms of functions defined on canonical intervals centered at the origin, we define the following mapped versions of $f$ and $g$
\begin{align*}
\tilde{f}(y) &= f\left(\frac{d-c}{2}y + \frac{a+b}{2}\right),\quad y \in [-(r+1),r+1], \\
\tilde{g}(y) &= g\left(\frac{d-c}{2}y + \frac{d+c}{2}\right),\quad y \in [-1,1],
\end{align*}
where $r = (b-a)/(d-c)-1 > 0$. Thus, the convolution of $\tilde{f}(y)$ and $\tilde{g}(y)$ is
\begin{align} 
\tilde{h}(y) = \int_{-1}^{1} \tilde{f}(y-t)\tilde{g}(t) \md t, \quad y \in [-r,r], \label{std}
\end{align}
which is related to $h(x)$ by
\begin{align*}
h(x) = \frac{d-c}{2} \tilde{h}\left(\frac{2x-(a+b+c+d)}{d-c}\right), \quad x\in[a+d,b+c].
\end{align*}
Thus, in the remainder of this paper we consider only the convolution $\tilde{h}$ given by \cref{std} with $\tilde{f}$ and $\tilde{g}$ defined on $[-(r+1),r+1]$ and $[-1,1]$. To further facilitate the discussion, we slightly abuse the notations by dropping the hats in $\tilde{f}$, $\tilde{g}$, and $\tilde{h}$, i.e., $f(x): [-(r+1), r+1] \rightarrow \mathbb{R}$, $g(x): [-1,1] \rightarrow \mathbb{R}$, and $h(x): [-r,r] \rightarrow \mathbb{R}$. Hence, we consider the convolution operator $\mF$ of Fredholm type that is defined by $f(x)$ and acts on $g(x)$
\begin{align*}
\mF_{f}[g](x) = h(x) = \int_{-1}^{1} f(x-t) g(t) \md t, \quad x \in [-r, r]. \label{fred}
\end{align*}
Here, $f(x)$ is usually called the kernel or the kernel function. Note that unlike in the convolution of Volterra type, the kernel function $f(x)$ and the function to be convolved, i.e., $g(x)$, do not commute. If $f(x)$ and $g(x)$ have certain regularity, they can be approximated by polynomials $f_M(x)$ and $g_N(x)$ of sufficiently high degree so that $\left\lVert f(x) - f_M(x) \right\rVert_{\infty}$ and $\left\lVert g(x) - g_N(x) \right\rVert_{\infty}$ are on the order of machine precision.
Such polynomial approximants can be constructed as, e.g., Legendre series
\begin{align}f_M(x) = \sum_{m=0}^M a_m P_m\left(\frac{x}{r+1}\right) ~~~~ \text{and} ~~~~ g_N(x) = \sum_{n=0}^N b_n P_n(x), \label{fM&gN}\end{align}
where $P_m$ is the Legendre polynomial of degree $m$. Suppose 
$h(x)$ is also approximated by a polynomial, denoted by $h_M(x)$. It follows then that $h_M(x)$ is the product of a $[-r,r] \times (N+1)$ quasi-matrix $\tilde{R}$ and the coefficient vector $\underline{b} = {\left(b_0, b_1 \ldots, b_N \right)}^{T}$. That is,
\begin{align}
h_M(x) =\tilde{R}\underline{b} = \left( 
      \vphantom{\begin{aligned}
          & \\ 
          & \\ 
          &
      \end{aligned}}
      \right.
      \begin{aligned}
          \mF_{f_M}[P_0]\Bigg| \mF_{f_M}[P_1] \Bigg| \cdots \Bigg| \mF_{f_M}[P_N]
      \end{aligned}
      \left.
      \vphantom{\begin{aligned}
          & \\ 
          & \\ 
          &
      \end{aligned}}
      \right)
   \underline{b}, \label{quasimat}
\end{align}
where the $n$th column of $\tilde{R}$ is the convolution of $f_M(x)$ and $P_n(x)$. For convenience, we denote by $\tilde{R}_n$ the $n$th column of $\tilde{R}$, with the index $n$ starting from $0$. Since $h_M(x)$ and $\tilde{R}_n$ are both polynomials of degree not exceeding $M$, they can be written as Legendre series, i.e.,
\begin{subequations}
\begin{align}
h_M(x) &= \sum_{m=0}^M c_m P_m\left(\frac{x}{r}\right) \label{hM}, \\
\tilde{R}_n(x) &= \sum_{m=0}^M R_{m,n} P_m\left(\frac{x}{r}\right), \label{Rn}
\end{align}\label{hM&Rn}
\end{subequations}
where $\{R_{m,n}\}_{m=0}^M$ are the Legendre coefficients of $\tilde{R}_n(x)$. Substituting \cref{hM&Rn} into \cref{quasimat} gives
\begin{align*}
\sum_{m=0}^{M} c_m P_m\left(\frac{x}{r}\right) = \sum_{n=0}^{N}b_n \sum_{m=0}^{M}R_{m,n}P_m\left(\frac{x}{r}\right),
\end{align*}
or, equivalently
\begin{align}
\underline{c} = R\underline{b},\label{cRb}
\end{align}
where $\underline{c} = {(c_0, \ldots, c_M)}^T$ and $R$ is an $(M+1)\times (N+1)$ matrix whose $(m,n)$th entry is $R_{m,n}$.

Since $f_M(x-t)$ is a polynomial in $t$ of degree at most $M$,
\begin{align}\tilde{R}_{n}(x)=\int_{-1}^{1}f_M(x-t)P_n(t)\,\md t=0 \label{fMPn}\end{align}
for $n>M$. Hence, we consider mainly the case of $M = N$ in \cref{sec:construction,sec:volterra}. In case of $N<M$, we prolong $\underline{b}$ by zeros. Equation \cref{fMPn} also justifies the employment of the Legendre polynomials in this problem.

At the time of writing, there is no literature establishing the construction of $R$, which is the primary goal of this work. In \cite{hal3}, Hale proposes the construction of $R$ via the approximation of two related convolution operators of Volterra type. However, Hale's method only works for the case $r=1$. 

The outline of this paper is as follows. In \cref{sec:construction}, we propose a novel approach to constructing the matrix approximation $R$ of the Fredholm convolution operator. 
\Cref{sec:volterra} relates the proposed approximation to that of the Volterra convolution operator. We demonstrate in \cref{sec:experiments} the speed and accuracy when the proposed approximations are employed to calculate convolution integrals, followed by various applications of such approximations. We close in \cref{sec:conclusion} with a summary.

\section{Constructing the convolution matrices}\label{sec:construction}
We start off by the lemma below which collects some basic results of Legendre polynomials.
\begin{lemma}\label{lemma1}
For $n\geq 1$,
\begin{subequations}
\begin{align}
(2n+1)P_n(x) &= \frac{\md}{\md x}\left(P_{n+1}(x)-P_{n-1}(x)\right), \label{diffP}\\
\int P_n(x) &= \frac{1}{2n+1}\left(P_{n+1}(x)-P_{n-1}(x)\right) +C, \label{intP} \\
(n+1)P_{n+1}(x) &= (2n+1)xP_n(x)-nP_{n-1}(x), \label{xP}
\end{align}\label{basic}%
\end{subequations}
where $C$ is the integration constant.
\end{lemma}
\begin{proof}
The proof can be found in many standard texts on Legendre polynomials, e.g., \cite{sze}.
\end{proof}

Given a Legendre series $S(x) = \displaystyle \sum_{l=0}^L \alpha_l P_l(x)$, it follows immediately from \cref{intP} that
\begin{align*}\int S(x) \md x = \sum_{l=0}^{L+1} \hat{\alpha}_l P_l(x) + C,\end{align*}
where $\left(\hat{\alpha}_0,\hat{\alpha}_1,\cdots\hat{\alpha}_{L+1}\right)^T = U_L \left(\alpha_0,\alpha_1,\cdots\alpha_{L}\right)^T$
and
\begin{align*}
U_L = \begin{pmatrix}
  0 & -\frac{1}{3}\\
  1 & 0 & -\frac{1}{5}\\
  & \frac{1}{3} & 0 & -\frac{1}{7}\\
  & & \ddots & \ddots &  \ddots 
\end{pmatrix}
\in\mathbb{R}^{(L+2)\times(L+1)}
\end{align*}
is the integration matrix. Similarly, it follows from \cref{xP} that the product of $x$ and $S(x)$
\begin{align*}x S(x) = \sum_{l=0}^{L+1} \tilde{\alpha}_l P_l(x),\end{align*}
where $\left(\tilde{\alpha}_0,\tilde{\alpha}_1,\cdots\tilde{\alpha}_{L+1}\right)^T = V_L \left(\alpha_0,\alpha_1,\cdots\alpha_{L}\right)^T$
and
\begin{align*}
V_L = \begin{pmatrix}
  0 &\frac{1}{3}\\
  1 &0 &\frac{2}{5}\\
  &\frac{2}{3} &0 &\frac{3}{7}\\
  & &\ddots &\ddots& \ddots 
\end{pmatrix}
\in\mathbb{R}^{(L+2)\times(L+1)}
\end{align*}
is the $x$-multiplication matrix. The following theorem shows that the matrix $R$ is skew upper triangular.
\begin{theorem} (skew upper triangularity) For $R$ given in \cref{cRb}, $R_{m,n} = 0$ for $m+n>M$.
\end{theorem}
\begin{proof}
Since $f_M$ is a polynomial of degree $M$, 
\begin{align*}f_M(x-t)=\sum_{l=0}^{M}Q_l(x)P_l(t),\end{align*}
where $\deg(Q_l)+\deg(P_l) \leq M$. For $m+n>M$,
\begin{align*}
R_{m,n} &= \frac{2m+1}{2r} \int_{-r}^{r}\int_{-1}^{1}f_M(x-t)P_n(t) \md t P_m\left(\frac{x}{r}\right) \md x\\
&= \frac{2m+1}{2r} \sum_{l=0}^{M}\int_{-r}^{r}P_m \left(\frac{x}{r}\right)Q_l(x) \md x \int_{-1}^{1}P_l(t)P_n(t) \md t = 0.
\end{align*}
To see why this must be zero, we assume that this is not the case. Then we must have $l = n$ and $\deg(Q_l) \geq m$ for some $l$. Thus, $M \geq \deg(Q_l) + l \geq m+n > M$, which is a contradiction.
\end{proof}

What follows next is our first main result which shows a four-term recurrence relation satisfied by the entries of $R$.
\begin{theorem}\label{thm:recur} (four-term recurrence relation)
For $m, n \geq 1$,
\begin{subequations}
\begin{align}
R_{m,n+1}&=R_{m,n-1}+r(2n+1)\left(\frac{R_{m-1,n}}{2m-1}-\frac{R_{m+1,n}}{2m+3} \right), \label{R_{m,n+1}}\\
R_{m+1,n}&=\frac{2m+3}{2m-1}R_{m-1,n}-\frac{1}{r}\frac{2m+3}{2n+1}\left(R_{m,n+1}-R_{m,n-1} \right), \label{R_{m+1,n}}\\
R_{m,n-1}&=R_{m,n+1}-r(2n+1)\left(\frac{R_{m-1,n}}{2m-1}-\frac{R_{m+1,n}}{2m+3} \right), \label{R_{m,n-1}}\\
R_{m-1,n}&=\frac{2m-1}{2m+3}R_{m+1,n}+\frac{1}{r}\frac{2m-1}{2n+1}\left(R_{m,n+1}-R_{m,n-1} \right). \label{R_{m-1,n}}
\end{align}\label{4term}
\end{subequations}
\end{theorem}
\begin{proof}
We show \cref{R_{m,n+1}} only, for \cref{R_{m+1,n},R_{m,n-1},R_{m-1,n}} are just paraphrases of \cref{R_{m,n+1}}. By a change of variable $z = y-t$, we have
\begin{align*}
\tilde{R}_{n+1}(y) = \int_{-1}^{1} f_M(y-t)P_{n+1}(t)\md t = \int_{y-1}^{y+1}f_M(z)P_{n+1}(y-z)\md z.
\end{align*}
Differentiating the last equation on both sides gives
\begin{align*}
  \frac{\md}{\md y}\tilde{R}_{n+1}(y) &= \int_{y-1}^{y+1}f_M(z)\frac{\md}{\md y}P_{n+1}(y-z)\md z + f_M(y+1)P_{n+1}(-1)\\
  & \hspace{5.5cm} -f_M(y-1)P_{n+1}(1)\\
  &= \int_{y-1}^{y+1}f_M(z)\left[(2n+1)P_n(y-z)+\frac{\md}{\md y}P_{n-1}(y-z)\right]\md z \\
  &\hspace{=2cm} + f_M(y+1)P_{n-1}(-1)-f_M(y-1)P_{n-1}(1)\\
  &= (2n+1)\tilde{R}_n(y)+\frac{\md}{\md y}\tilde{R}_{n-1}(y),
\end{align*}
where we have used \cref{diffP} and $P_{n+1}(-1)=P_{n-1}(-1)$  and $P_{n+1}(1)=P_{n-1}(1)$. Integrating on both sides over $[-r, x]$ yields
\begin{align*}\tilde{R}_{n+1}(x)=(2n+1)\int_{-r}^{x}\tilde{R}_n(y)\md y+\tilde{R}_{n-1}(x)+C.\end{align*}
Substituting \cref{Rn} into the last equation and noting the triangularity of $R$, we have
\begin{align*}
\sum_{m=0}^{M-n-1}R_{m,n+1}P_m\left(\frac{x}{r}\right)=&(2n+1) \sum_{m=0}^{M-n}R_{m,n}\int_{-r}^{x}P_m\left(\frac{y}{r}\right)\md y\\
    &+ \sum_{m=0}^{M-n+1}R_{m,n-1}P_m\left(\frac{x}{r}\right)+C.
\end{align*}
Matching the Legendre coefficients of $P_m\left(\frac{x}{r}\right)$ for $m\geq 1$ leads to \cref{R_{m,n+1}}.
\end{proof}

The four identities in \cref{4term} suggest that the entries in $R$ can be calculated recursively in four ways: columnwise rightward, rowwise downward, columnwise leftward, and rowwise upward. It may not be obvious to immediately see how we can start up the recursion with \cref{R_{m,n-1}} or \cref{R_{m-1,n}} which we leave for the moment. To recurse with \cref{R_{m,n+1}} or \cref{R_{m+1,n}}, it is apparent that we have to first figure out the first two columns and the zeroth row or the first two rows and the zeroth column respectively. To this end, we first give the following lemma on representing $P_j\left(\frac{x+1}{r+1}\right)-P_j\left(\frac{x-1}{r+1}\right)$ in terms of $\left\{P_{k}\left(\frac{x}{r}\right)\right\}_{k=0}^{j-1}$. For notational convenience, we let $\underline{P}_j = \left( P_0\left(\frac{x}{r}\right), P_1\left(\frac{x}{r}\right), \cdots, P_j\left(\frac{x}{r}\right) \right)$.

\begin{lemma}\label{lem:col}
Suppose that 
\begin{align*}
P_j\left(\frac{x+1}{r+1}\right)-P_j\left(\frac{x-1}{r+1}\right)=\sum_{k=0}^{j-1}w_{k,j}P_{k}\left(\frac{x}{r}\right).
\end{align*}
If the coefficients $w_{k,j}$ are collected in the matrix $W_L=(w_{k,j})\in\mathbb{R}^{(L+1)\times(L+2)}$, $W_L$ is strictly upper triangular with
\begin{align*}
W_L(\texttt{0:2,0:3})=
\begin{pmatrix}
0&\frac{2}{r+1}&0&\frac{2r^2-6r+2}{{(r+1)}^3}\\
0&0&\frac{6r}{{(r+1)}^2}&0\\
0&0&0&\frac{10r^2}{{(r+1)}^3}\\
\end{pmatrix}.
\end{align*}
With the further assumption that $w_{k,j} = 0$ for $k<0$ for all $j$, the entries of $W_L$ satisfy a $9$-term recurrence relation
\begin{equation}
  \begin{aligned}
    w_{k,j+2}= B _{k-1,j+1}w_{k-1,j+1}+C_{k+1,j+1}w_{k+1,j+1}+E_{k-2,j}w_{k-2,j}+F_{k,j}w_{k,j}\\
    +G_{k+2,j}w_{k+2,j}+H_{k-1,j-1}w_{k-1,j-1}+J_{k+1,j-1}w_{k+1,j-1}+K _{k,j-2}w_{k,j-2},\\
  \end{aligned}
  \label{9term}
\end{equation}
where
\begin{gather*}
A_{k,j+2} = \frac{2j+3}{(j+2)(r+1)},~B_{k-1,j+1}=\frac{2rkA_{k,j+2}}{2k-1},~ C_{k+1,j+1}=\frac{2r(k+1)A_{k,j+2}}{2k+3},\\
E_{k-2,j}=-\frac{r^2k(k-1)(2j+1)A_{k,j+2}}{(r+1)(2k-3)(2k-1)(j+1)},~ F_{k,j} =F_{k,j}^1+F_{k,j}^2+F_{k,j}^3,\\
F_{k,j}^1 = \frac{(2j+1)A_{k,j+2}}{(j+1)(r+1)},~F_{k,j}^2 = -\frac{(r+1)(j+1)A_{k,j+2}}{2j+3}-\frac{(r+1)j^2A_{k,j+2}}{(j+1)(2j-1)},\\
F_{k,j}^3 = -\frac{r^2(2j+1)(2k^2-2k-1)A_{k,j+2}}{(r+1)(j+1)(4k^2-4k-3)},\\
G_{k+2,j} = -\frac{r^2(2j+1)(k+1)(k+2)A_{k,j+2}}{(r+1)(j+1)(2k+5)(2k+3)},~H_{k-1,j-1}=\frac{2rjkA_{k,j+2}}{(j+1)(2k-1)},\\
J_{k+1,j-1} = \frac{2rj(k+1)A_{k,j+2}}{(j+1)(2k+3)},~K_{k,j-2} = -\frac{j(j-1)(r+1)A_{k,j+2}}{(j+1)(2j-1)}.
\end{gather*}
Additionally, $w_{k,j}=0$ when $k+j$ is even.
\end{lemma}
\begin{proof}
Let $\phi_j(x)=P_j\left(\frac{x+1}{r+1}\right)-P_j\left(\frac{x-1}{r+1}\right)$ and $\psi_j(x)=P_j\left(\frac{x+1}{r+1}\right)+P_j\left(\frac{x-1}{r+1}\right)$. It is easy to see that $\phi_j(x)$ is a polynomial of degree $j-1$ for $j\geq 1$, implying that $W_L$ is strictly upper triangular. The first few $\phi_j(x)$ can be easily obtained with a bit of algebraic manipulation:
\begin{gather*}
\phi_0(x) = 0,\quad \phi_1(x) = \frac{2}{r+1},\quad \phi_2(x) = \frac{6r}{{(r+1)}^2} P_1\left(\frac{x}{r}\right),\\
\phi_3(x) = \frac{2r^2-6r+2}{{(r+1)}^3}+ \frac{10r^2}{{(r+1)}^3}P_2\left(\frac{x}{r}\right).
\end{gather*}
It follows from \cref{xP} that
\begin{align*}
\phi_{j+1}(x)=&\frac{2j+1}{j+1}\left[\frac{x+1}{r+1}P_{j}\left(\frac{x+1}{r+1}\right)-\frac{x-1}{r+1}P_{j}\left(\frac{x-1}{r+1}\right)\right]\\
&-\frac{j}{j+1}\left[P_{j-1}\left(\frac{x+1}{r+1}\right) - P_{j-1}\left(\frac{x-1}{r+1}\right)\right]\\
&= \frac{2j+1}{(j+1)(r+1)}\left(x\phi_j(x) + \psi_j(x)\right) - \frac{j}{j+1}\phi_{j-1}(x).
\end{align*}
Similarly,
\begin{align*}
\psi_{j+1}(x)=\frac{2j+1}{(j+1)(r+1)}\left( x\psi_j(x) + \phi_j(x)\right) - \frac{j}{j+1}\psi_{j-1}(x).
\end{align*}
Combining the last two equations, we have
\begin{align}
\phi_{j+2}(x)&=\frac{2j+3}{(j+2)(r+1)}\left[2x\phi_{j+1}(x)+\frac{2jx}{j+1}\phi_{j-1}(x)-\frac{j(j-1)(r+1)}{(j+1)(2j-1)}\phi_{j-2}(x)\right.\nonumber \\
&+\left.\left(\frac{2j+1}{(j+1)(r+1)}\left(1-x^2\right)-\frac{(j+1)(r+1)}{2j+3}-\frac{j^2(r+1)}{(j+1)(2j-1)}\right)\phi_j(x)\right], \label{phirec}
\end{align}
which relates $\phi_{j-2}(x), \ldots, \phi_{j+2}(x)$ together. To get rid of the multiplying factors of $x$ and $x^2$, we note
\begin{equation}
\begin{aligned}
\left(\frac{x}{r}\right)\phi_j(x) = \underline{P}_j V_{j-1}{\left(w_{0,j},w_{1,j}\cdots w_{j-1,j}\right)}^T,\\
{\left(\frac{x}{r}\right)}^2 \phi_j(x) = \underline{P}_{j+1} V_{j}V_{j-1}{\left(w_{0,j},w_{1,j}\cdots w_{j-1,j}\right)}^T. \label{x/r}
\end{aligned}
\end{equation}
Substituting \cref{x/r} into \cref{phirec} and matching up the Legendre coefficients in the resulting equation gives \cref{9term}. In case of $k+j$ being even, $w_{k,j} = 0$ follows from the $9$-term recurrence relation.
\end{proof}

We are now in a position to construct the zeroth column of $R$.
\begin{theorem} \label{thm:0th_col}(Construction of the zeroth column) The elements of the zeroth column of $R$ are
\begin{align}(r+1)W_M U_M \underline{a}, \label{0thcol}\end{align}
where $\underline{a} = {\left(a_0, a_1, \ldots, a_M\right)}^T$.
\end{theorem}
\begin{proof}
The zeroth column of $R$ is the vector collecting the coefficients of $\int_{-1}^{1}f_M(x-t)\md t$ in the basis ${\{P_m\left(\frac{x}{r}\right)\}}_{m=0}^M$. Since it is a polynomial in $x$ of degree at most $M$, 
\begin{align}
\int_{-1}^{1}f_M(x-t)\md t &= \int_{-1}^{1} \sum_{m=0}^{M} a_m P_m\left(\frac{x-t}{r+1}\right)\md t \notag\\
&=(r+1)\sum_{m=0}^{M+1} \hat{a}_m \left[P_m\left(\frac{x+1}{r+1}\right)-P_m\left(\frac{x-1}{r+1}\right)\right],\label{subeq}
\end{align}
where $\underline{\hat{a}}={\left(\hat{a}_0,\cdots \hat{a}_{M+1}\right)}^T = U_M \underline{a}$.
It follows from \cref{lem:col} that
\begin{align*}
\sum_{m=0}^{M+1}\hat{a}_m\left[P_m\left(\frac{x+1}{r+1}\right)-P_m\left(\frac{x-1}{r+1}\right)\right] = \underline{P}_M W_M  \underline{\hat{a}},
\end{align*}
which, combined with \cref{subeq}, gives
\begin{align*}
\int_{-1}^{1}f_M(x-t)\,d t = (r+1)\underline{P}_M W_M U_M \underline{a}.
\end{align*}
\end{proof}

The following theorem shows the construction of the first column of $R$. 
\begin{theorem}\label{thm:1st_col} (construction of the first column) The elements of the first column of $R$ are
\begin{align}
\left(r (r+1) V_{M+1}W_M U_M - (r+1)^2 W_{M+1} U_{M+1} V_M\right)\underline{a}. \label{1stcol}
\end{align}
\end{theorem}

\begin{proof}
Since the entries in the first column of the matrix $R$ are the coefficients of $\int_{-1}^{1}f_M(x-t)t \md t$,
\begin{align*}
&\int_{-1}^{1}f_M(x-t)t \md t
= (r+1) \sum_{m=0}^{M}a_m \int_{-1}^{1} P_m \left(\frac{x-t}{r+1}\right)\left(\frac{x}{r+1} - \frac{x-t}{r+1}\right) \md t\\
=& x\sum_{m=0}^{M}a_m \int_{-1}^{1} P_m\left(\frac{x-t}{r+1}\right)\md t -(r+1) \sum_{m=0}^{M}a_m \int_{-1}^{1} \frac{x-t}{r+1} P_m\left(\frac{x-t}{r+1}\right)\md t. \label{1st2terms}
\end{align*}

By \cref{thm:0th_col} and \cref{xP},
\begin{align*}
\sum_{m=0}^{M}a_m \int_{-1}^{1}\frac{x-t}{r+1} P_m\left(\frac{x-t}{r+1}\right)\md t =(r+1)\underline{P}_{M+1} W_{M+1} U_{M+1}V_M \underline{a}.
\end{align*}
For the first term, up to a constant $r$,
\begin{align*}
&\frac{x}{r}\sum_{m=0}^{M}a_m \int_{-1}^{1} P_m\left(\frac{x-t}{r+1}\right)\md t = (r+1)\underline{P}_{M+1} V_{M+1}W_M U_M \underline{a},
\end{align*}
which follows from an analogous argument. The last three equations gives \cref{1stcol}.
\end{proof}
The length of \cref{1stcol} is $M+2$, whose last two entries are exactly zero, since $\tilde{R}_{1}(x)$ is a polynomial of degree $M-1$. 

The zeroth and first rows can be constructed similarly. \cref{lem:hatw} and \cref{thm:0th_row} below are the results in parallel for rows, and we omit the proofs.

\begin{lemma}\label{lem:hatw}
Suppose that
\begin{align*}
P_j\left(\frac{r-t}{r+1}\right)-P_j\left(-\frac{r+t}{r+1}\right) = \sum_{k=0}^{j-1}\hat{w}_{k,j}P_k(t),
\end{align*}
If the coefficients $\hat{w}_{k,j}$ are collected in the matrix $\hat{W}_L=(\hat{w}_{k,j})\in\mathbb{R}^{(L+1)\times(L+2)}$, $\hat{W}_L$ is strictly upper triangular with
\begin{align*}
\hat{W}_L(\texttt{0:2,0:3})=
\begin{pmatrix}
0&\frac{2r}{r+1}&0&\frac{r(2r^2-6r+2)}{{(r+1)}^3}\\
0&0&-\frac{6r}{{(r+1)}^2}&0\\
0&0&0&\frac{10r}{{(r+1)}^3}
\end{pmatrix},    
\end{align*}
With the further assumption that $\hat{w}_{k,j} = 0$ for $k<0$ for all $j$, the entries of $\hat{W}_L$ satisfy a $9$-term recurrence relation
\begin{equation*}
\begin{aligned}
  \hat{w}_{k,j+2}= \hat{B}_{k-1,j+1}\hat{w}_{k-1,j+1}+\hat{C}_{k+1,j+1}\hat{w}_{k+1,j+1}+\hat{E}_{k-2,j}\hat{w}_{k-2,j}+\hat{F}_{k,j}\hat{w}_{k,j}\\
  +\hat{G}_{k+2,j}\hat{w}_{k+2,j}+\hat{H}_{k-1,j-1}\hat{w}_{k-1,j-1}+\hat{J}_{k+1,j-1}\hat{w}_{k+1,j-1}+\hat{K}_{k,j-2}\hat{w}_{k,j-2},\\
\end{aligned}
\end{equation*}
where

\begin{gather*}
\hat{B}_{k-1,j+1}=-B_{k-1,j+1}/r,~ \hat{C}_{k+1,j+1}=-C_{k+1,j+1}/r,~ \hat{E}_{k-2,j}=E_{k-2,j}/r^2,\\
\hat{F}_{k,j} = \hat{F}_{k,j}^1+\hat{F}_{k,j}^2+\hat{F}_{k,j}^3,~ \hat{F}_{k,j}^1 = r^2F_{k,j}^1,~\hat{F}_{k,j}^2 = F_{k,j}^2,~\hat{F}_{k,j}^3 = F_{k,j}^3/r^2,\\
\hat{G}_{k+2,j}=G_{k+2,j}/r^2,~ \hat{H}_{k-1,j-1}=-H_{k-1,j-1}/r,\\
\hat{J}_{k+1,j-1}=-J_{k+1,j-1}/r,~\hat{K}_{k,j-2}=K_{k,j-2}.
\end{gather*}
Additionally, $\hat{w}_{k,j} = 0$ when $k+j$ is even.
\end{lemma}

For the first two rows, we also need the diagonal matrix $D_M = \diag\left(1,1/3,\right.$ $\left.\ldots,1/(2M+1)\right) \in \mathbb{R}^{(M+1)\times (M+1)}$ which comes from the integration constant. 
\begin{theorem}\label{thm:0th_row}
The elements of the zeroth and the first rows of $R$ are
\begin{align}
&\frac{r+1}{r}D_M \hat{W}_M U_M \underline{a}, \label{0throw} \\
&\frac{3(r+1)}{r^2}D_{M+1}\left((r+1)\hat{W}_{M+1}U_{M+1}V_M + V_{M+1} \hat{W}_M U_M\right)\underline{a}, \label{1strow}
\end{align}
respectively.
\end{theorem}

\begin{figure}[t!]
\centering
\includegraphics[width=0.5\linewidth]{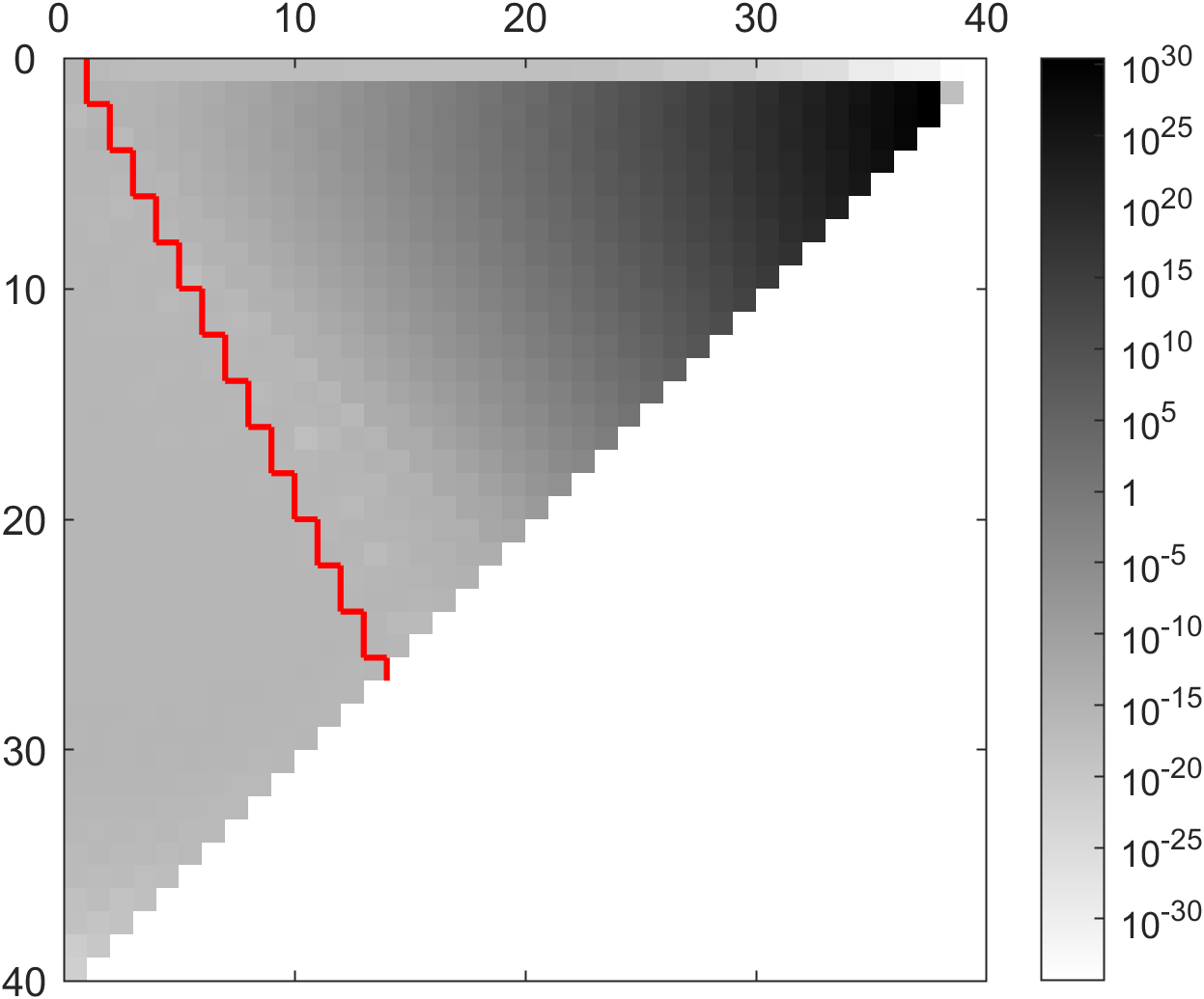}
\caption{Error growth when recursing for $R$ corresponding to $f_{39}$ using \cref{R_{m,n+1}}, which is stable in the region that is to the left of the red borderline.}\label{fig:unstable}
\end{figure}

The calculation of $R$ could have been as easy as suggested by \cref{R_{m,n+1}} and \cref{R_{m+1,n}}---generate the first two columns and the first row or the first two rows and the first column followed by recursing the entries columnwise or rowwise respectively. Unfortunately, \cref{R_{m,n+1},R_{m+1,n}} are numerically unstable. In fact, recursing with \cref{R_{m,n+1}} and \cref{R_{m+1,n}} only works stably when $\frac{r(2n+1)}{2m-1} \leq 1$ and $\frac{2m+3}{r(2n+1)} \leq 1$ respectively. Otherwise, the error in $R_{m\pm 1,n}$ and $R_{m,n\pm 1}$  are to be amplified repeatedly in the course of recursion. In the worst-case scenario, the error in certain entries could be amplified approximately by a factor of $M\mathbb{!}$. To see the instability, consider $R \in \mathbb{R}^{40\times 40}$ defined by $f_{39}(x)$ for $r=2$, which is a Legendre series of degree $39$ with all the coefficients set to $1$. \cref{fig:unstable} shows the entrywise error in $R$ obtained by recursing columnwise with \cref{R_{m,n+1}}. The largest error is $2.7098 \times 10^{30}$ in this example, occurring at $R_{1, 37}$.

\begin{figure}[t!]
\centering
\subfloat[$r \geq 1$]{\includegraphics[width=0.3\linewidth]{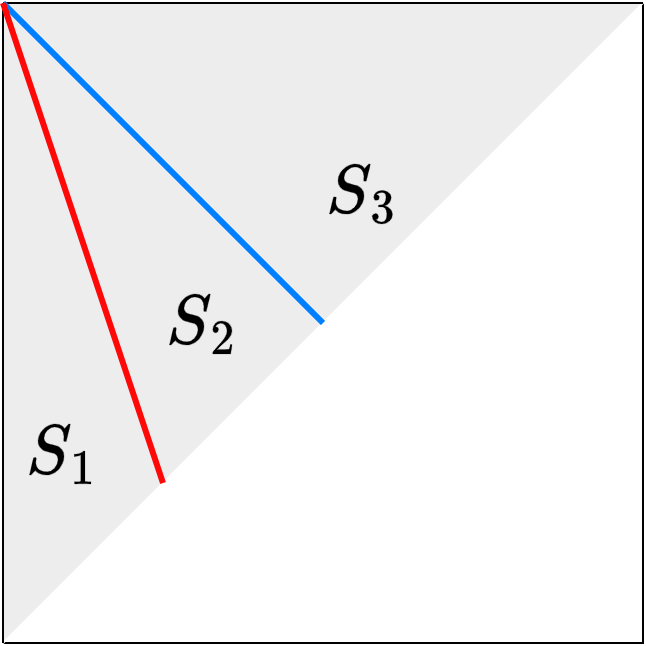}\label{fig:r>1}}
\hspace{2cm}
\subfloat[$0 < r \leq 1$]{\includegraphics[width=0.3\linewidth]{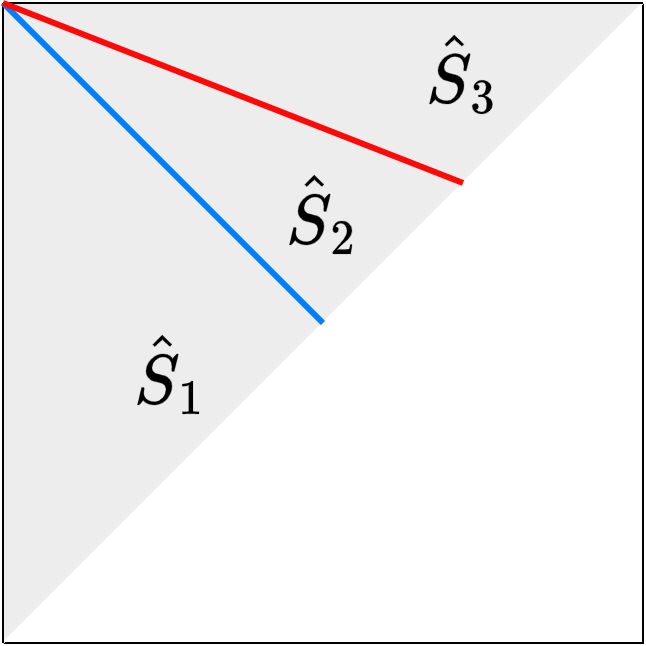}\label{fig:r<1}}
\caption{Partition of the nonzero region by the borderline of stability (red) and the line of dependence (blue).}\label{fig:partition}
\end{figure}

To identify the entries for which each of \cref{4term} is stable, we first safely ignored the amplification effect caused by $\frac{2m+3}{2m-1}$ in \cref{R_{m+1,n}} since it is close enough to $1$ for large $m$. Thus, the four recurrence relations in \cref{4term} are stable respectively when
\begin{align}
r\frac{2n+1}{2m-1}\leq 1,~~~ \frac{2m+3}{r(2n+1)}\leq 1,~~~ r\frac{2n+1}{2m-1}\leq 1,~~~ \frac{2m-1}{r(2n+1)}\leq 1. \label{stable}
\end{align}
If we further ignore the constants in the fractions in \cref{stable}, \cref{R_{m,n+1},R_{m,n-1}} are stable when $m/n \geq r$ whereas \cref{R_{m-1,n},R_{m+1,n}} are stable for $m/n \leq r$. The stability regions for \cref{4term} are marked by red lines in \cref{fig:partition}. For $r\geq 1$, \cref{R_{m,n+1},R_{m,n-1}} are stable in $S_1$, whereas \cref{R_{m+1,n},R_{m-1,n}} are stable in $S_{23} = S_2+S_3$. For $0 < r \leq 1$, \cref{R_{m,n+1},R_{m,n-1}} are stable in $\hat{S}_{12} = \hat{S}_1+\hat{S}_2$, while \cref{R_{m+1,n},R_{m-1,n}} are stable in $\hat{S}_3$. 

\begin{figure}[t!]
\centering
\subfloat[$r \geq 1$]{\includegraphics[width=0.45\linewidth]{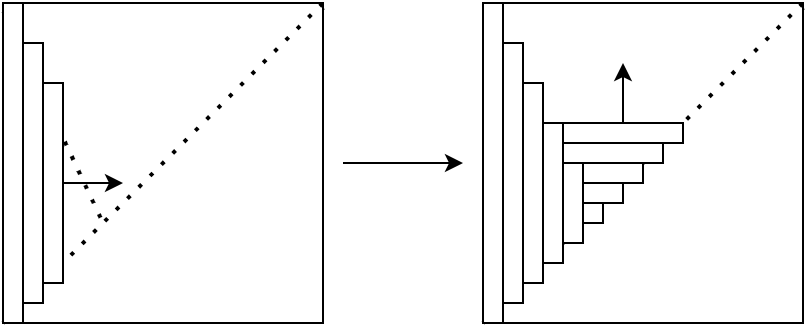}\label{fig:prog1}}
\hfill
\subfloat[$0< r \leq 1$]{\includegraphics[width=0.45\linewidth]{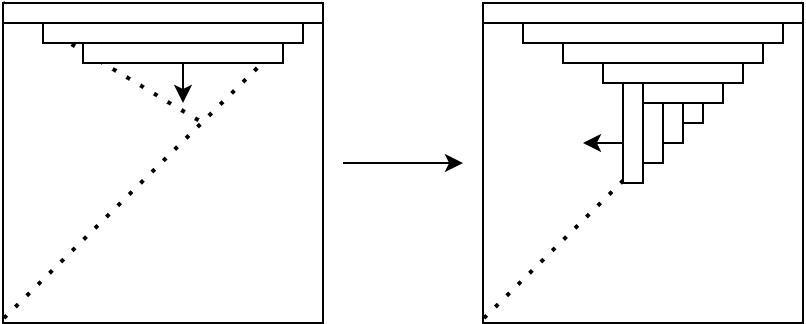}\label{fig:prog2}}
\caption{Stable recursing strategies.}\label{fig:prog}
\end{figure}

Now it should be clear that we have four possible recursing strategies for entries of $R$ when $r\geq 1$: (1) recursing with \cref{R_{m,n+1}} rightward in $S_1$ followed by recursing with \cref{R_{m-1,n}} upward in $S_{23}$, (2) recursing with \cref{R_{m+1,n}} downward in $S_{23}$ followed by recursing with  \cref{R_{m,n-1}} leftward in $S_1$, (3) recursing with \cref{R_{m,n+1}} rightward in $S_1$ followed by recursing with \cref{R_{m+1,n}} downward in $S_{23}$, and (4) recursing with \cref{R_{m+1,n}} downward in $S_{23}$ followed by recursing with \cref{R_{m,n+1}} rightward in $S_1$. We would not prefer (3) and (4) since we have to precompute both the first two columns and rows, which incurs unnecessary cost. If we take (2), we have to know certain entries in $S_1$ to recurse downward for the entries in $S_2$. Hence, when $r>1$ we recurse following strategy (1) (\cref{fig:prog1}). Following an analogous argument, we recurse with \cref{R_{m+1,n}} downward in $\hat{S}_3$ followed by recursing with \cref{R_{m,n-1}} leftward in $\hat{S}_{12}$ for $0 \leq r \leq 1$ (\cref{fig:prog2}). We summarize the stable recursions in \cref{alg:stable}. If we recalculate the example in \cref{fig:unstable} using \cref{alg:stable}, the largest entrywise error is about $2.3 \times 10^{-16}$.

\begin{algorithm}[t!]
\caption{Stable construction of the convolution matrix $R$}\label{alg:stable}
\begin{algorithmic}
\IF{$r\geq 1$}
  \STATE Construct the first two columns by \cref{0thcol,1stcol}. \hfill $\rhd 42M^2$
  \STATE{Recursing with \cref{R_{m,n+1}} rightward for the entries in $S_1$.} \hfill $\rhd 6M^2/(r+1)$
  \STATE{Recursing with \cref{R_{m-1,n}} upward for the entries in $S_{23}$.} \hfill $\rhd 7rM^2/(r+1)$
\ELSE
  \STATE{Construct the first two rows by \cref{0throw,1strow}.} \hfill $\rhd 42M^2$
  \STATE{Recursing with \cref{R_{m+1,n}} downward for the entries in $\hat{S}_3$.} \hfill $\rhd 7rM^2/(r+1)$
  \STATE{Recursing with \cref{R_{m,n-1}} leftward for the entries in $\hat{S}_{12}$.} \hfill $\rhd 6M^2/(r+1)$
\ENDIF
\end{algorithmic}
\end{algorithm}

Since each nonzero entry of $W_M$ or $\hat{W}_M$ is obtained at a cost of approximately $145$ flops and only $M^2/4$ entries of $W_M$ or $\hat{W}_M$ are nonzero, the total complexity of constructing $W_M$ or $\hat{W}_M$ is about $36 M^2$. If we ignore the $\mO(M)$ operations, the dominant cost in evaluating \cref{0thcol,1stcol} or \cref{0throw,1strow} is the three multiplication of $W_M$ or $\hat{W}_M$ with vectors, totaling $6M^2$ flops. The cost for generating the first two columns or rows is thus $42M^2$. For the $M^2/(2(r+1))$ entries in $S_1$ or $\hat{S}_{12}$, recursing for each entry with \cref{R_{m,n+1}} or \cref{R_{m,n-1}} requires $12$ flops, leading to a total cost of $6M^2/(r+1)$ flops. Analogously, each entry in $S_{23}$ or $\hat{S}_{3}$ can be obtained following \cref{R_{m-1,n}} and \cref{R_{m+1,n}} respectively with $15$ flops. For the $rM^2/(2(r+1))$ entries in $S_{23}$ or $\hat{S}_{3}$, the total cost is $7rM^2/(r+1)$. Thus, the total complexity in the recursing part of \cref{alg:stable} is roughly $7M^2$, independent of $r$. We list the stepwise complexities in \cref{alg:stable} and note that the overall cost of \cref{alg:stable} is roughly $49M^2$ flops.

\section{Relation to the convolution operator of Volterra type}\label{sec:volterra}
For $\check{f}, \check{g}:[-1,1]\rightarrow \mathbb{R}$, the convolution operator of Volterra type is defined as
\begin{align*}
\check{h} = \mV_{\check{f}}[\check{g}](x) = \int_{-1}^x \check{f}(x-t)\check{g}(t) \md t, \label{Vf}
 \end{align*}
where the resulting function $\check{h}:[-2,0]\rightarrow \mathbb{R}$. Suppose that $\check{f}(x)$ and $\check{g}(x)$ are approximated by the Legendre series $\check{f}_M(x) = \sum_{m = 0}^M \check{a}_m P_m(x)$ and $\check{g}_N(x) = \sum_{n = 0}^N \check{b}_n P_n(x)$ respectively. Then $\check{h}_{M+N+1}(x) = \mF_{\check{f}_M}[\check{g}_N](x)$ is a Legendre series of length $M+N+2$, and the $(M+N+2)\times (N+1)$ matrix that relates $\check{h}_{M+N+1}(x)$ and $\check{g}_N(x)$ serves as an approximation of $\mV_{\check{f}}$. If we denote by $V_{\check{f}}$ this matrix approximation of $\mV_{\check{f}}$, $V_{\check{f}}$ is also obtainable by a four-term recurrence relation similar to \cref{4term} at a cost of $\mO((M+N)N)$. Particularly, $V_{\check{f}}$ is banded with the upper and lower bandwidths both being exact $M+1$. For detail, see \cite{hal2,xu2}.

For a function $\varphi$, let $\varphi_{\{a,b\}}$ denote the restriction of $\varphi$ to the interval $[a,b]$ and $\hat{\varphi} = \varphi(-x)$. It is shown in \cite{hal3} that when $r=1$ a convolution operator of Fredholm type can be written as the combination of two Volterra convolution operators
\begin{align}
\mF_f[g](x) = \mV_{f_1} [g](x) + \mV_{\hat{f}_2}[\hat{g}](\hat{x}), \label{FVV}
\end{align}
where $f_1 = f_{\{0,2\}}$ and $f_2 = f_{\{-2,0\}}$. Let $\hat{I} = \text{diag}(1, -1, 1, -1, \cdots)$ be the alternating diagonal matrix. Suppose that both $f_1$ and $\hat{f}_2$ can be approximated by Legendre series of degree $M$ and $g$ by a Legendre series of degree $N$. Hence, the $(M+N+2)\times (N+1)$ matrix approximation of the Fredholm convolution operator 
\begin{align}
\hat{R} = V_{f_1} + \hat{I}V_{\hat{f}_2}\hat{I}. \label{RRV}
\end{align}
Note that $\hat{R}$ may have different dimensions from $R$, as the latter is assumed square throughout. Obtaining $\hat{R}$ via \cref{RRV} incurs the construction of both $V_f$ and $V_{\hat{f}}$, the cost being $\mO((M+N)N)$. On the contrary, the cost of constructing $R$ following \cref{alg:stable} is $\mO(M^2)$, independent of $N$. 

When $\hat{R}$ is constructed via \cref{RRV}, cancellations are expected for entries in the band except those in the top-left triangular part. Since the dimensions of $V_{f}$ and $V_{\hat{f}}$ are usually determined by the length of $g_N(x)$, a great deal of the computation might be wasteful when $N$ is much larger than $M$. Furthermore, the cancellation is never perfect---entrywise residuals are usually on the order of machine epsilon or greater. When such a $\hat{R}$ is used to calculate convolutions or solve convolution integral equations, the cancellation errors may lead to inaccurate, if not totally erroneous, results. This is particularly the case if certain entries in the nonzero triangular part are on the same order of or smaller than the cancellation errors. 

\begin{figure}[t]
\centering
\subfloat[$\hat{R}$ obtained by \cref{RRV}]{\includegraphics[width=0.44\linewidth]{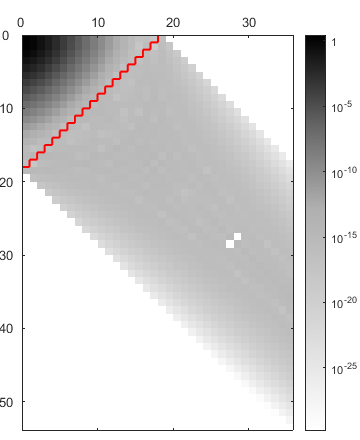}\label{fig:sum}}
\hfill 
\subfloat[$R$ obtained by \cref{alg:stable}]{\includegraphics[width=0.44\linewidth]{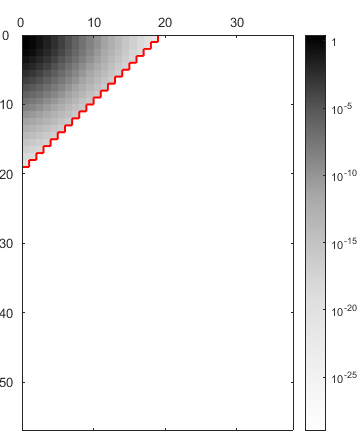}\label{fig:rec}}
\caption{Entrywise magnitude of $\hat{R}$ and $R$.}\label{fig:plus}
\end{figure}

Let $f_{17}(x)$ be the Legendre series of degree $17$ which approximates $e^x$ on $[-2,2]$ and consider the Fredholm convolution operator defined by $f_{17}(x)$. \cref{fig:sum} shows the entrywise magnitude of $\hat{R}$ obtained by \cref{RRV}, where the nonzero entries that could be obtained directly by \cref{alg:stable} are enclosed by the borderline in red. Comparing with \cref{fig:rec} which displays the entrywise magnitude of $R$ constructed using \cref{alg:stable}, we can see that the cancellation errors are of the same magnitude as or even larger than many nonzero entries, making the skew upper-triangular structure much invisible. This may explain why the triangularity is overlooked by Hale and Townsend in their works \cite{hal3,hal2}.

For $r\neq 1$, there is unfortunately no direct relation like \cref{FVV}. Hence, \cref{alg:stable} is so far the only way to construct approximations to the Fredholm convolution operator (see \cref{sec:conveqn}). However, the Fredholm convolution of two functions can still be calculated via convolutions of Volterra type \cite{hal2}. For $r < 1$, it is shown that 
\begin{align}
\mF_{f}[g](x) = \mV_{f_{\{r-1,r+1\}}}[g](x) + \mV_{\hat{f}_{\{-r-1,r-1\}}}\left[\hat{g}_{\{-2r+1,1\}}\right](\hat{x}),~~~x \in  [-r,r]. \label{r<1}
\end{align}
Thus, we need to make two restrictions of $f$ and one of $g$. These are done by evaluating $f$ and $g$ at Chebyshev grids of $M$ and $N$ points respectively via the Clenshaw algorithm at a cost of $\mO(M^2+N^2)$, followed by Chebyshev--Legendre transforms to recover the Legendre coefficients for each restriction in $\mO(M(\log M)^2 + N(\log N)^2)$ flops. Finally, constructing two Volterra convolution matrices costs $\mO\left((M+N)M\right)$ flops.

For $r > 1$, 
\begin{align}
\mF_{f}[g](x) =\sum_{j=1}^{\left\lfloor r\right\rfloor } \mF_{f^{[j]}}[g](x) + \mF_{f_{\{2\left \lfloor r\right\rfloor -r-3, r+1\} }}[g](x),~~~x \in [-r,r], \label{r>1}
\end{align}
where $f^{[j]} = f_{\{-(r+1)+2(j-1), -(r+1)+2(j+1)\}}$ and the summation and the plus sign should be understood as concatenations. Each summand $\mF_{f^{[j]}}[g](x)$ satisfies $r=1$ since the length of the domain of $f^{[j]}$ is twice that of $g$. When $r$ is an integer, the last term in \cref{r>1} vanishes. In such a case, we need to make $r$ restrictions of $f$ and construct $2r$ approximations of the Volterra convolution operators, incurring $\mO\left(rM(M+N+(\log M)^2)\right)$ flops. If $r$ is not an integer, the last term $\mF_{f_{\{2\left \lfloor r\right\rfloor -r-1, r+1\}}}[g](x)$ corresponds to $r<1$ and can be treated following \cref{r<1}, adding an extra $\mO(N^2 + N (\log N)^2)$ term to the total complexity.

In contrast to the complexities we just worked out, the cost of \cref{alg:stable} 
is independent of $N$ and $r$. This immediately suggests the advantage of the proposed approximation in calculating the convolution for which the domain of $f$ is much sizable than that of $g$ or $N$ is much greater than $M$. See \cref{sec:computingconv} for numerical experiments.

\section{Numerical examples}\label{sec:experiments}
In this section, we first test the use of the approximation of the Fredholm convolution operator for calculating the convolution of two given functions (\cref{sec:computingconv,sec:weierstrass}), followed by its application in solving convolution integral equations (\cref{sec:conveqn}) 
and the calculation of pseudospectra (\cref{sec:convps}). All the numerical experiments are performed in \textsc{Julia} v1.9.3 on a laptop with a 14 core 4.7 Ghz Intel Core i7 CPU and 16GB RAM. The execution times are obtained using \texttt{BenchmarkTools.jl}.

\subsection{Speed in computing Fredholm convolutions}\label{sec:computingconv}
\begin{figure}[t!]
\centering
\subfloat[$r=1$]{\includegraphics[width=0.48\linewidth]{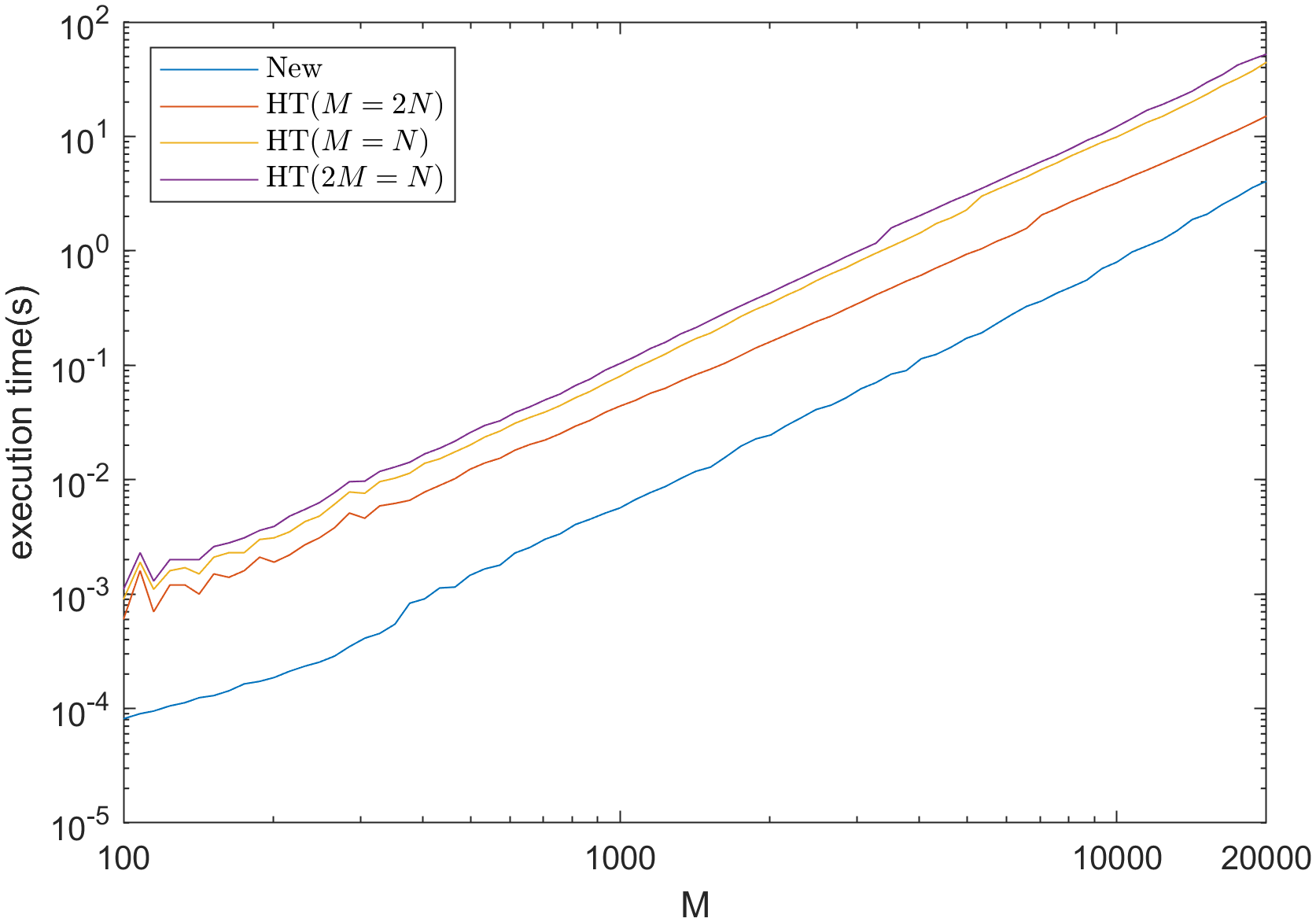}}
\hfill 
\subfloat[$r=2$]{\includegraphics[width=0.48\linewidth]{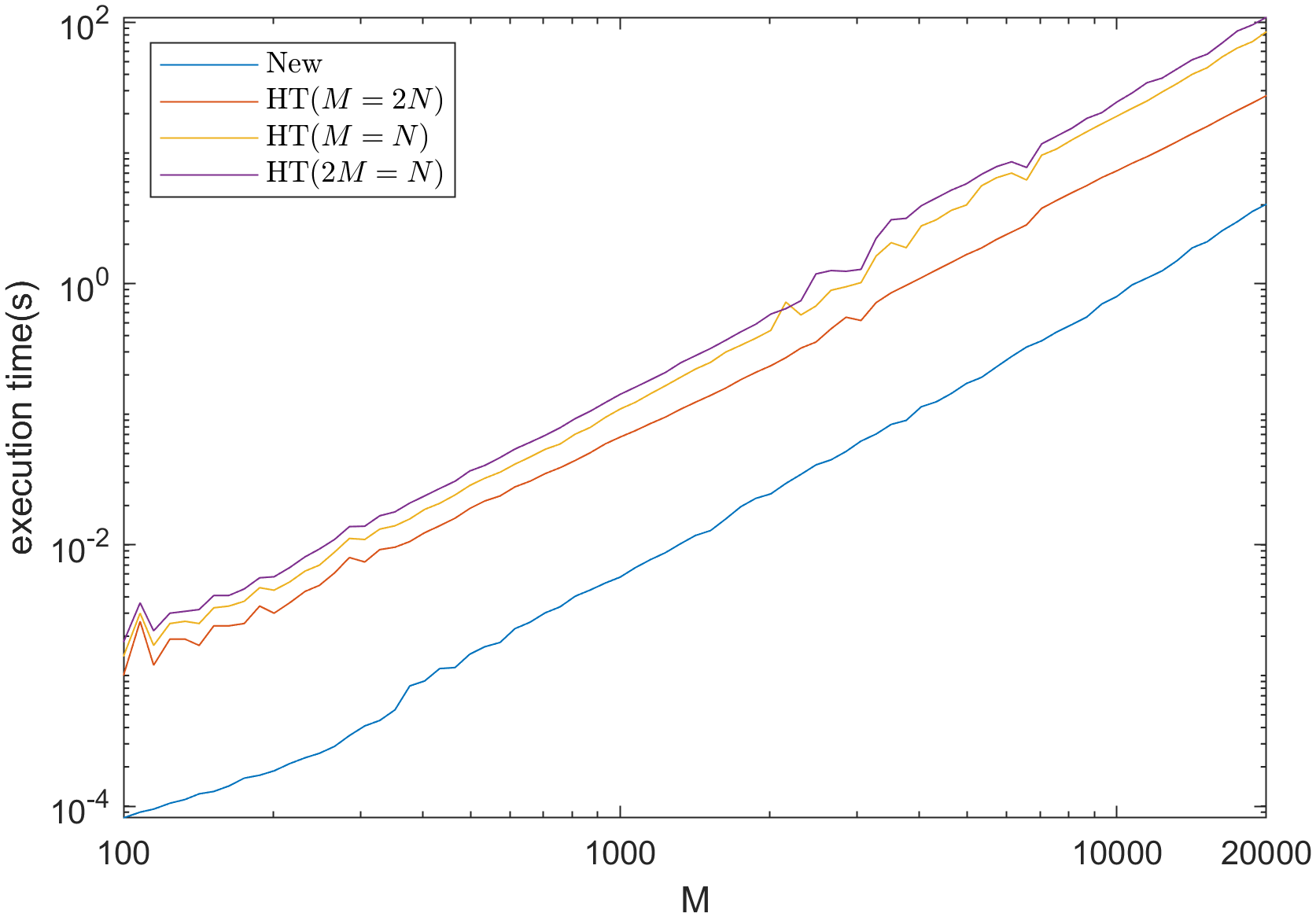}}
\\
\subfloat[$r=10$]{\includegraphics[width=0.48\linewidth]{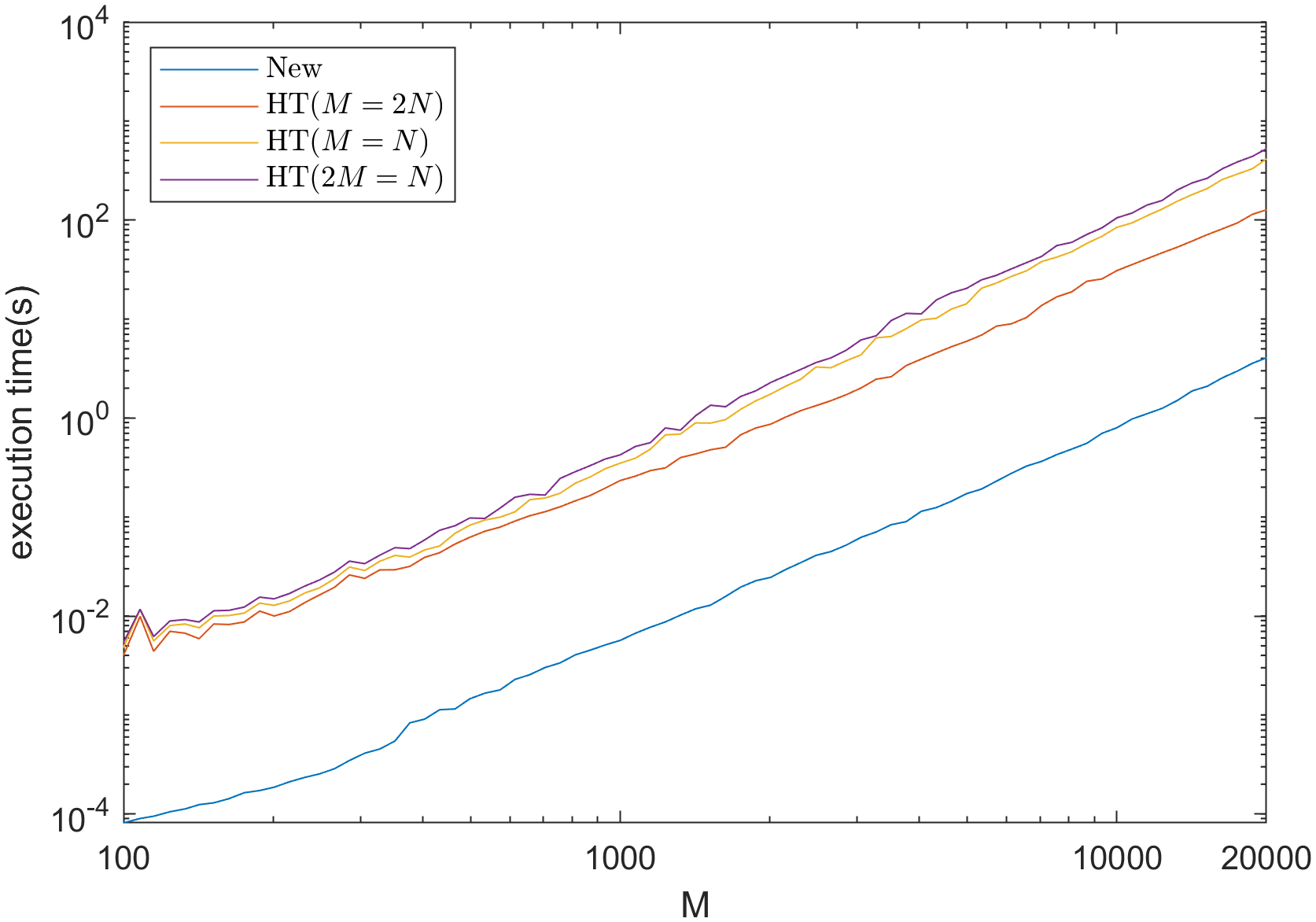}}
\hfill 
\subfloat[$r=100$]{\includegraphics[width=0.48\linewidth]{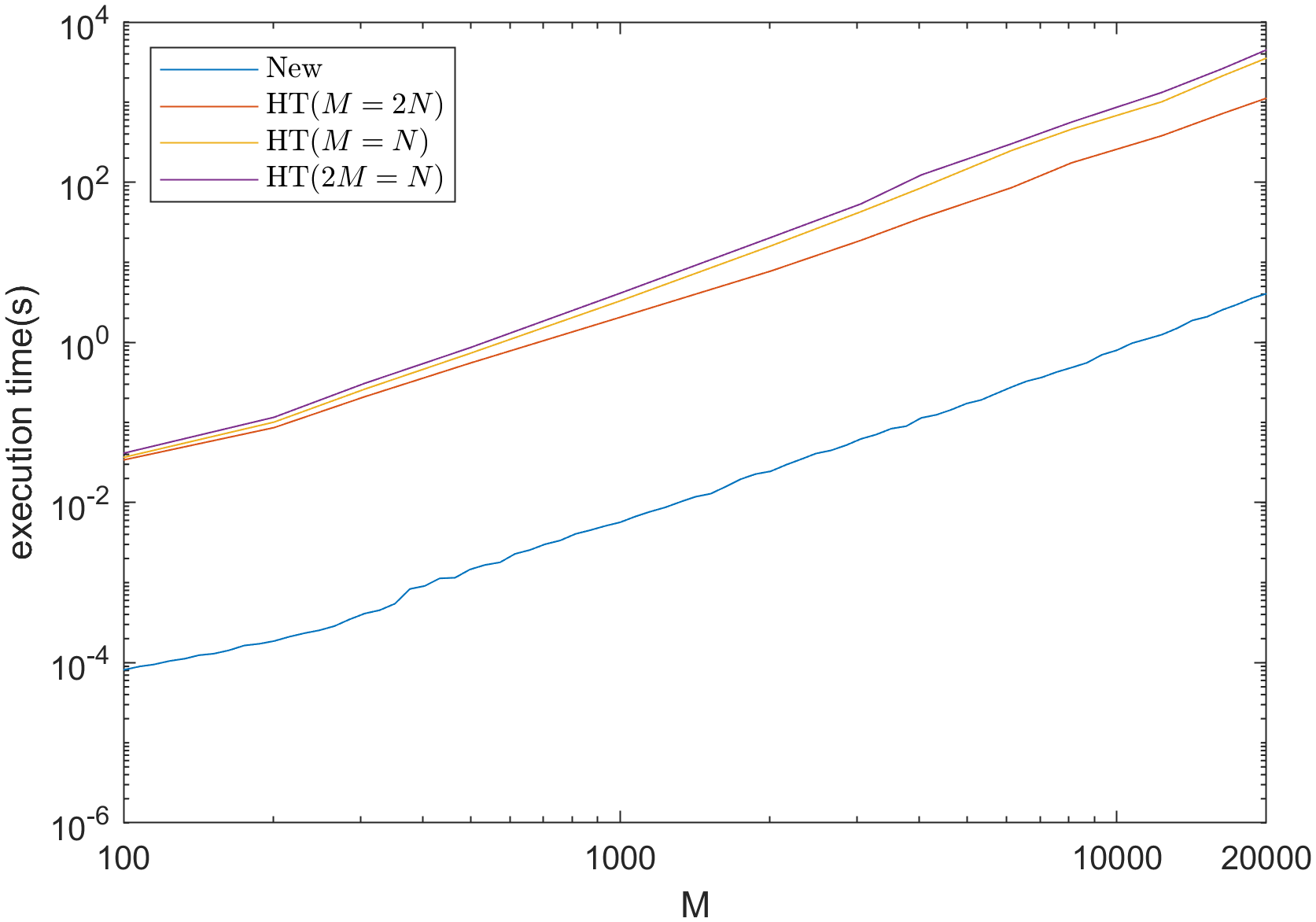}}
\caption{Computational time versus $M$ for various choices of $N$ and $r$.}\label{fig:M}
\end{figure}

We evaluate the Fredholm convolution of two given Legendre series with randomly generated coefficients using the proposed approximation and compare the computational time with that of the Hale--Townsend method \cite[\S 5]{hal2}. Since the Hale--Townsend method depends on $M$, $N$, and $r$, we time the computation by varying one of them and fixing the other two. 

Each panel of \cref{fig:M} has fixed value of $r$ and shows the dependence of the execution time on $M$ for various values of $N/M$. It is shown that the time grows quadratically with $M$ for both methods, as expected. The dependence on $N$ by the Hale--Townsend method is suggested by the curves corresponding different ratios of $N/M$, whereas the proposed method is $N$-independent. Given that the proposed method is independent of $r$, its execution times remain consistent and are thus indistinguishable across the panels. On the other hand, the time for the Hale--Townsend method grows proportionally with $r$. In the case of $r = 1$, the new method based on the proposed approximation is $5$ to $10$ times faster than the Hale--Townsend method. For $r=100$, the speed-up is roughly $200$ folds. 

\begin{figure}[t!]
\centering
\subfloat[$r=1$]{\includegraphics[width=0.48\linewidth]{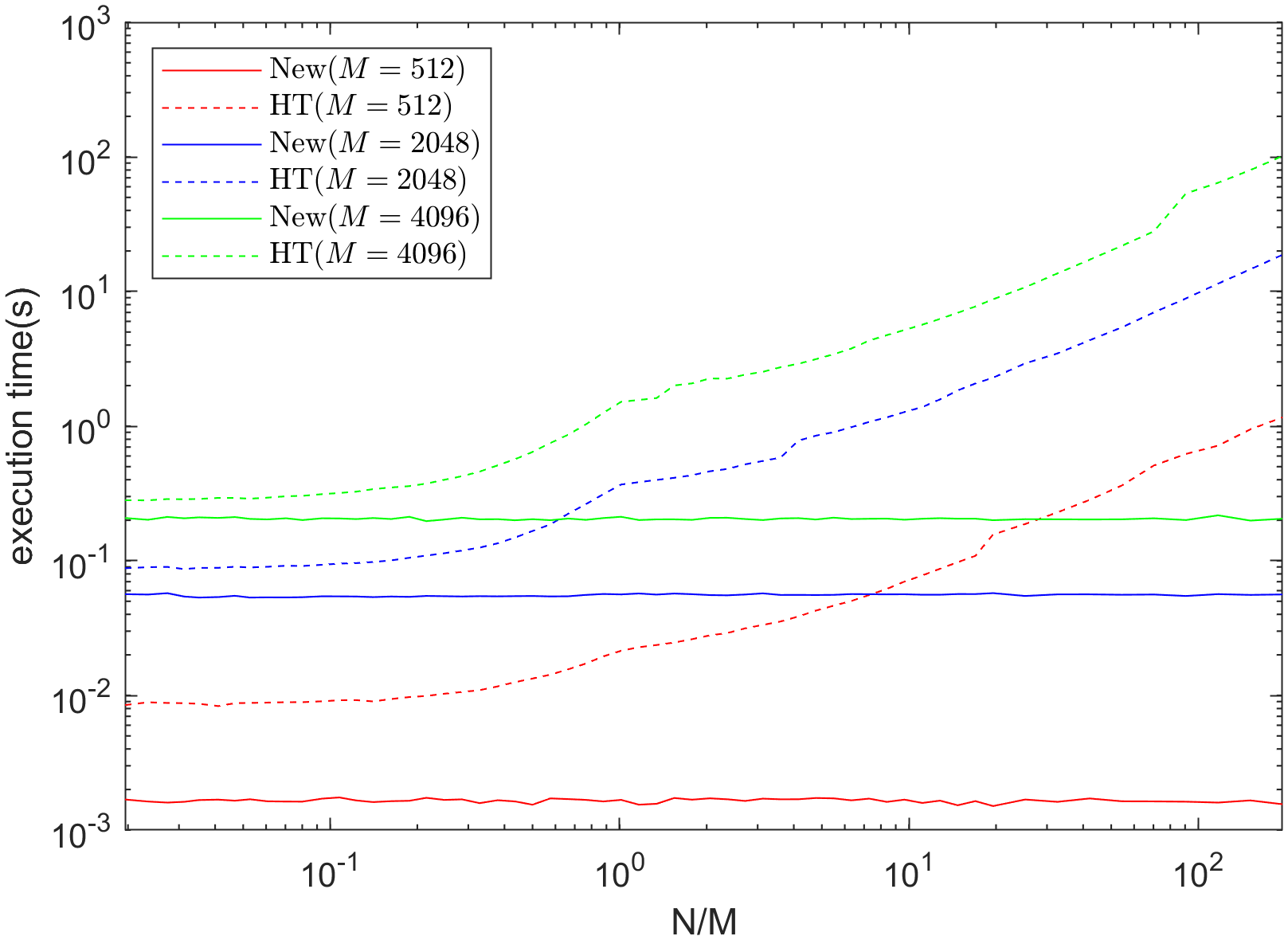}}
\hfill 
\subfloat[$r=2$]{\includegraphics[width=0.48\linewidth]{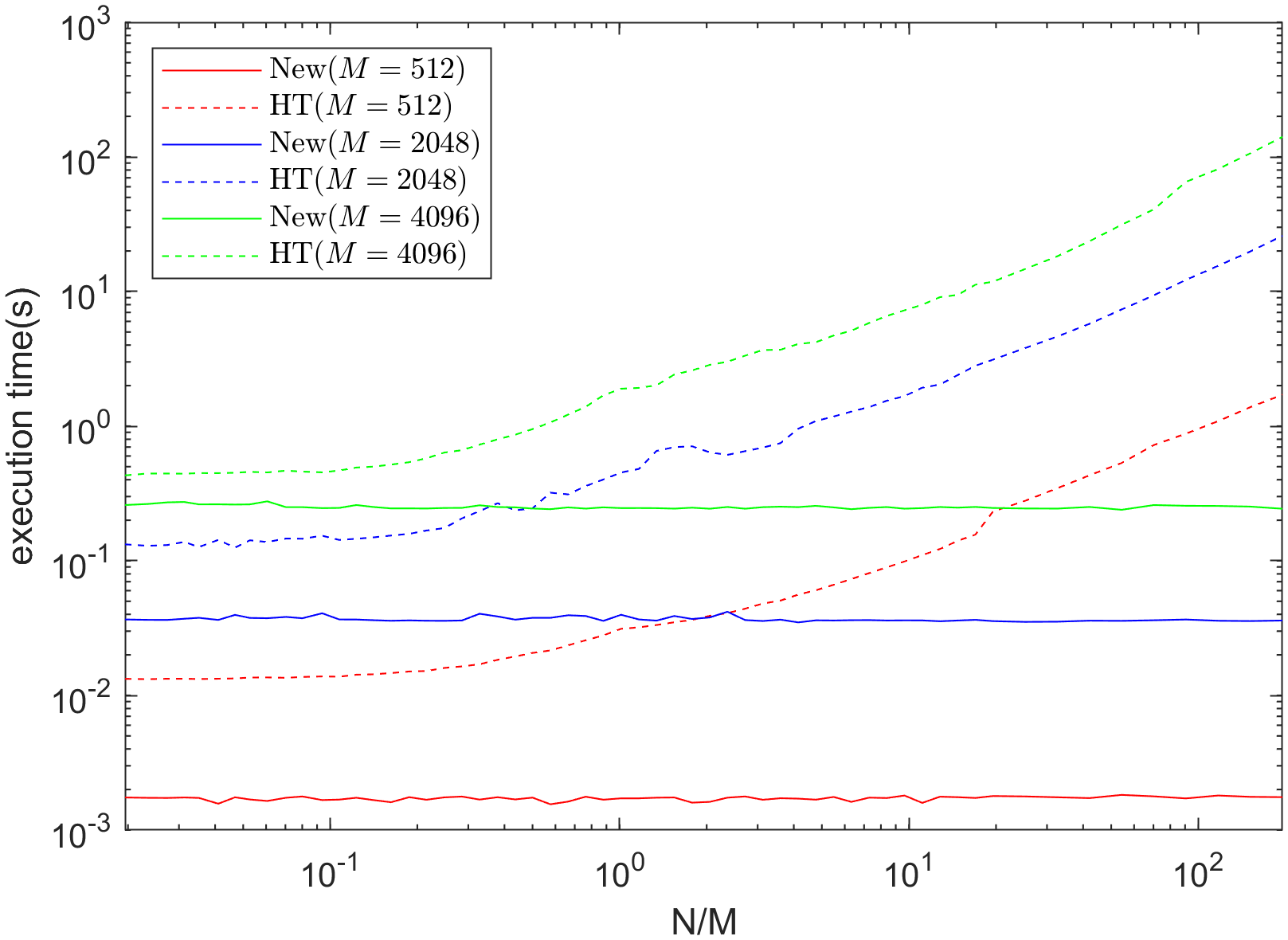}}
\\
\subfloat[$r=10$]{\includegraphics[width=0.48\linewidth]{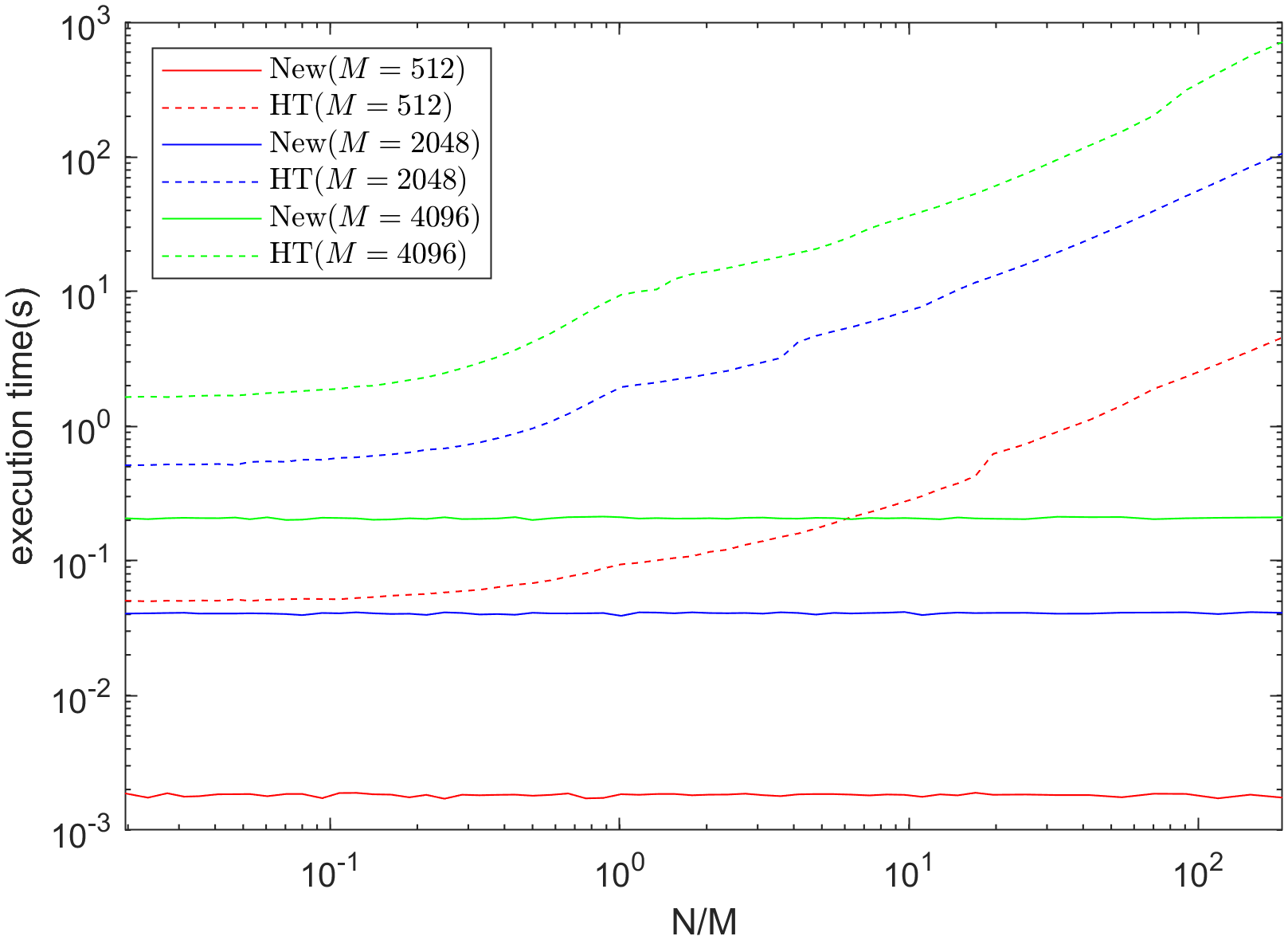}}
\hfill 
\subfloat[$r=100$]{\includegraphics[width=0.48\linewidth]{n_r10.png}}
\caption{Computational time versus $N/M$ for various choices of $M$ and $r$.}\label{fig:N/M}
\end{figure}

\begin{figure}[t!]
\centering
\subfloat[$M=10N$]{\includegraphics[width=0.48\linewidth]{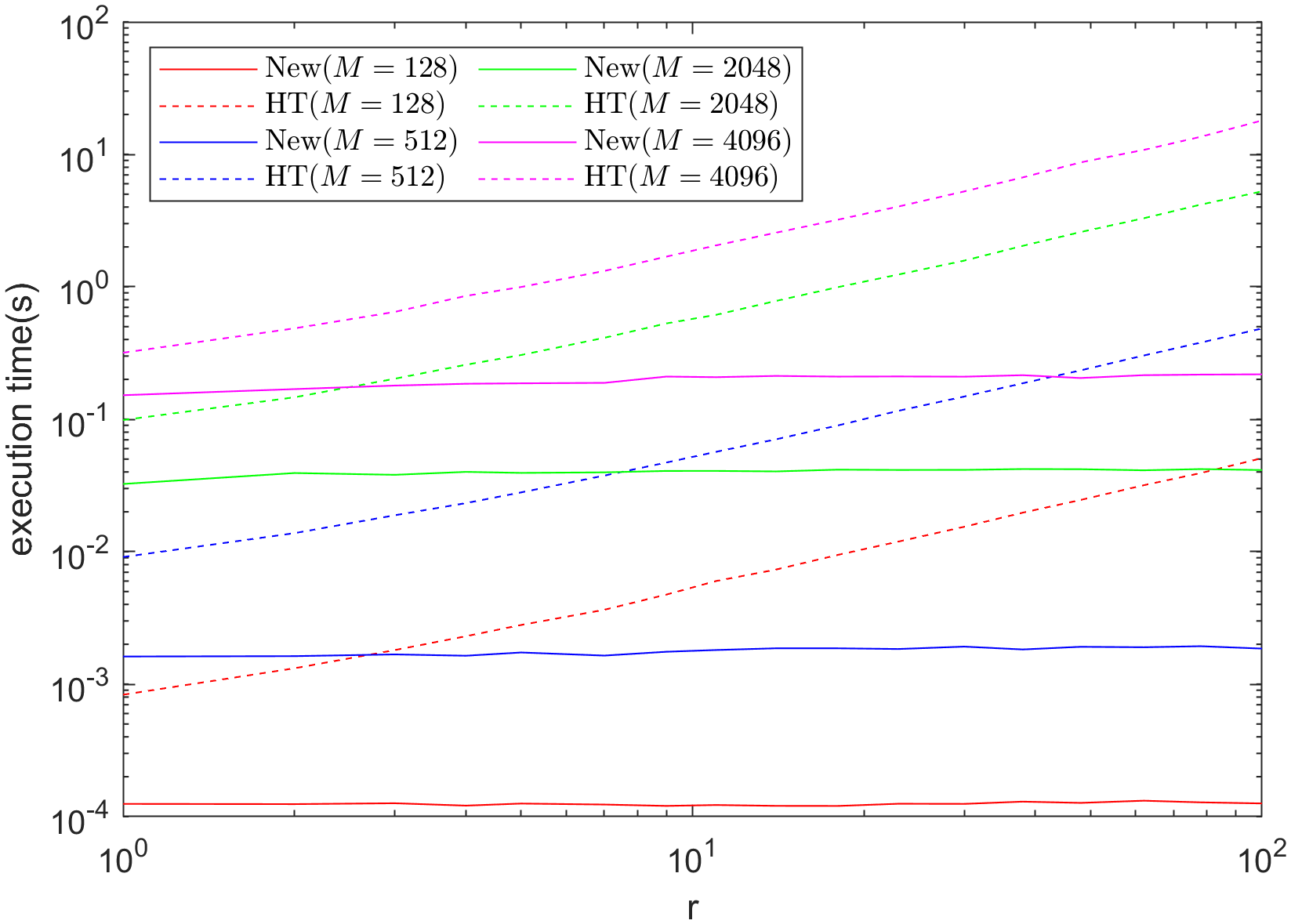}}
\hfill 
\subfloat[$M=2N$]{\includegraphics[width=0.48\linewidth]{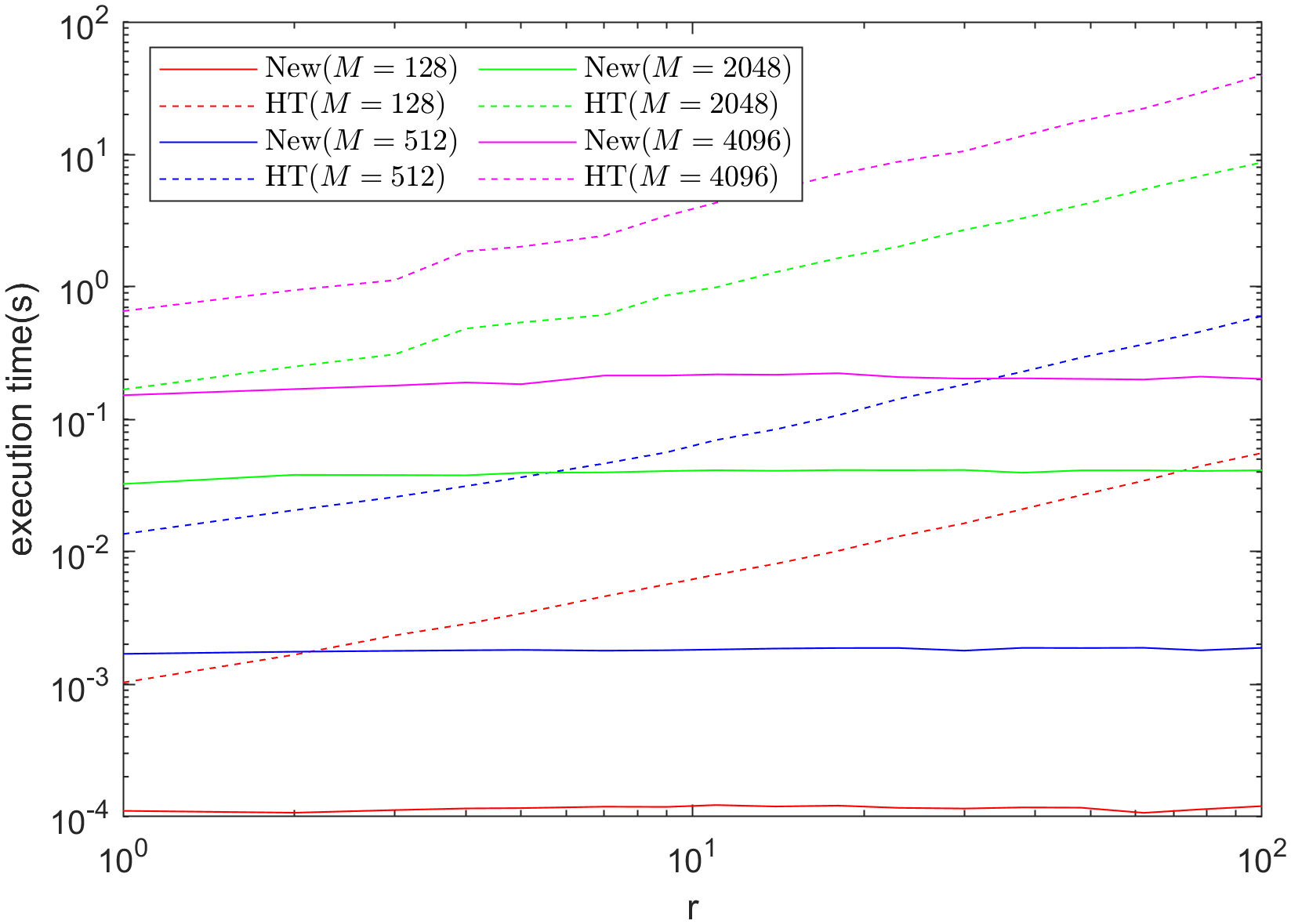}}
\\
\subfloat[$M=N$]{\includegraphics[width=0.48\linewidth]{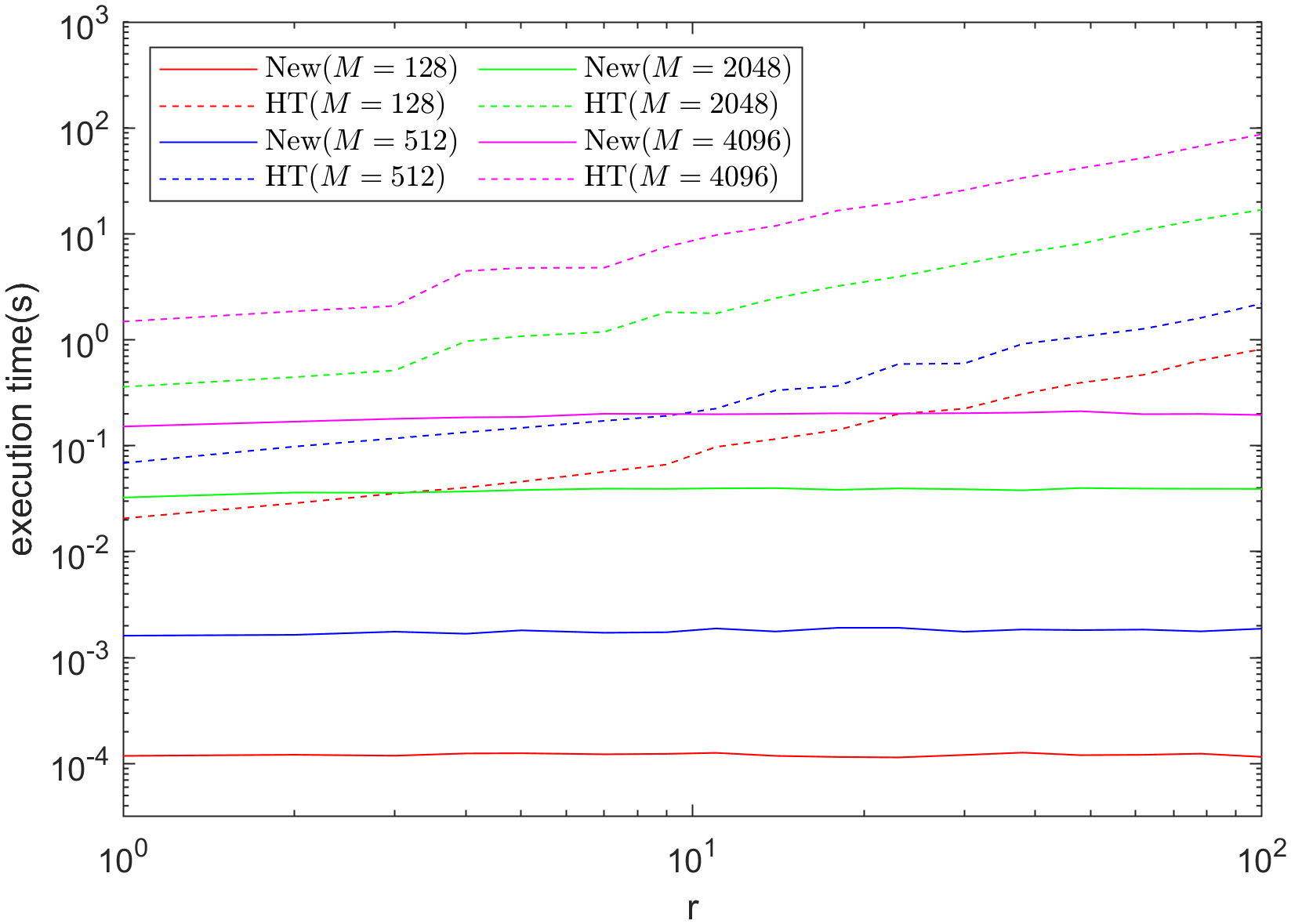}}
\hfill 
\subfloat[$10M=N$]{\includegraphics[width=0.48\linewidth]{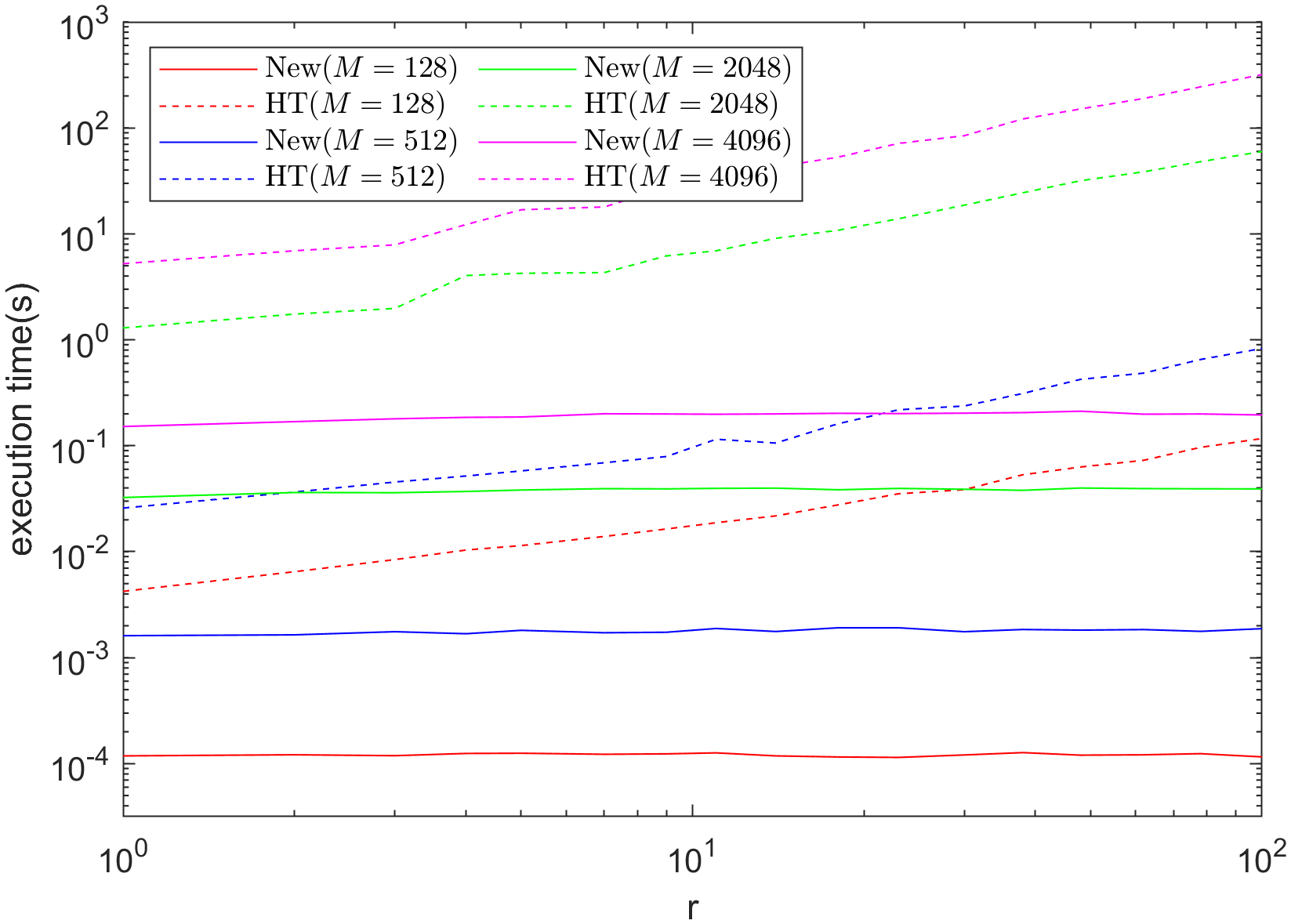}}
\caption{Computational time versus $r$ for various choices of $M$ and $N/M$.}\label{fig:r}
\end{figure}

\cref{fig:N/M} shows the time for varying $N/M$ ratio at fixed $M$ and $r$. We could have shown the plots for time versus $N$. However, the use of the ratio $N/M$ helps align the curves for different $M$ so that the dependence of $M$ is better illustrated. Since the proposed method is $N$-independent the solid curves are largely straight lines, whereas the time for the Hale--Townsend method scales linearly with $N$. The greater $N$ is, the more conspicuous the speed-ups are gained by the new method. \cref{fig:r} displays the (in)dependence on $r$ of the two methods. This is again where the proposed method is favored as its computational time is constant despite the length of the intervals. In \cref{fig:M,fig:N/M,fig:r}, $M$, $N$, and $r$ are chosen to be typical and representative, and the new method is faster in all these tests with clear advantages for $r \gg 1$ or $N \gg M$.

\subsection{An example for demonstrating the Weierstrass approximation theorem}\label{sec:weierstrass}
\begin{figure}[t!]
\centering
\subfloat[Weierstrass]{\includegraphics[width=0.48\linewidth]{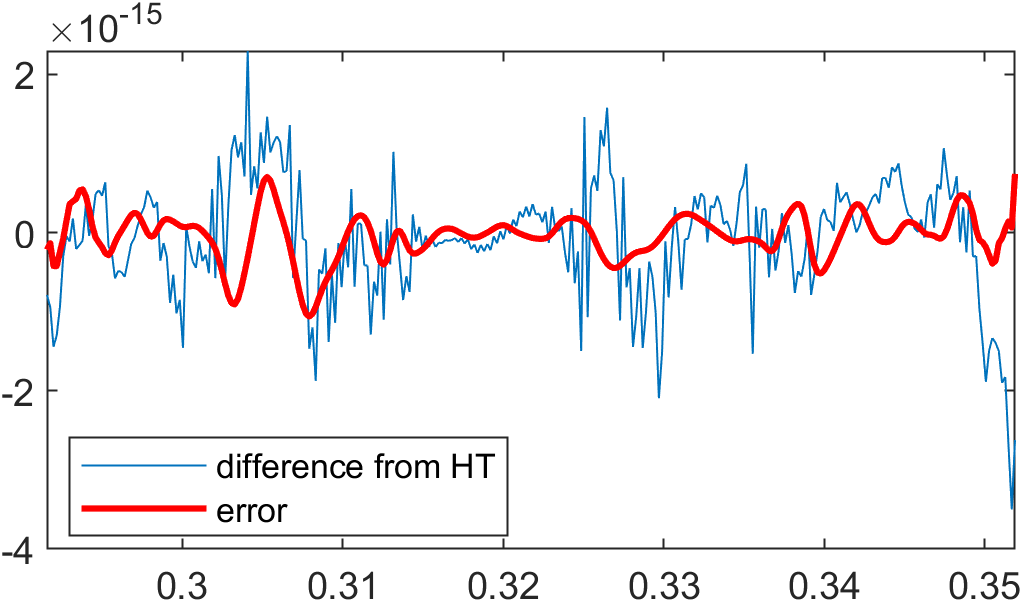} \label{fig:atap}}
\hfill 
\subfloat[Love's equation]{\includegraphics[width=0.48\linewidth]{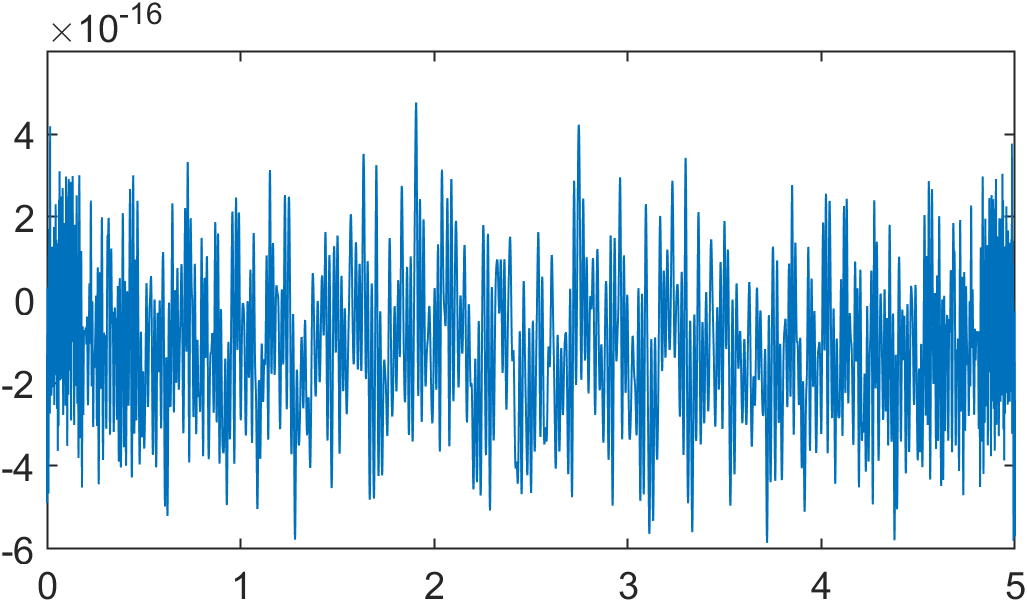}\label{fig:love}}
\caption{(a) The difference of the computed convolution from the result of the Hale--Townsend method and the result obtained using extended precision. (b) Error in the computed solution to Love's equation.}\label{fig:wl}
\end{figure}

Our second example is the evaluation of the Fredholm convolution of $f(x) = \sin(1/x)\sin(1/(\sin(1/x)))$ with $x \in [0.2885554757,0.3549060246]$ and $g(x) = \exp(-x^2/(4t))/\sqrt{4\pi t}$ with $t=10^{-7}$ and $x\in [-3\times 10^{-3}, 3\times 10^{-3}]$. This example is borrowed from \cite[\S 6]{ATAP}, where it is used to demonstrate the Weierstrass approximation theorem. As in \cite[\S 6]{ATAP}, we approximate $f(x)$ and $g(x)$ by polynomials of degree $2000$ and $65$ respectively. \cref{fig:atap} shows the error obtained by comparing with the result from computing in octuple precision\footnote{The $256$-bit octuple precision is the default format of \textsc{Julia}'s \texttt{BigFloat} type of floating point number. \texttt{BigFloat}, based on the GNU MPFR library, is the arbitrary precision floating point number type in \textsc{Julia}.} and the difference from the result obtained using the Hale--Townsend method. The new method is more than twice faster than the Hale--Townsend method in this example though $N$ is much smaller than $M$. 

\subsection{Convolution integral equation}\label{sec:conveqn}
The approximation of Fredholm convolution operators is the key in solving convolution integral equations of Fredholm type. For instance, Love's integral equation \cite{lov}
\begin{align}
y(t) = f(t) - \frac{\delta}{\pi} \int_0^1 \frac{y(s)}{\delta^2+{(t-s)}^2} \md s, ~~~~ 0\leq t\leq 5 \label{love}
\end{align}
is a Fredholm convolution integral equation of the second kind that arises in electrostatics. Here, we choose $\delta = -1$ and $f(t) = 1 +\frac{\delta}{\pi}\left(\arctan(1-t)+\arctan(t)\right)$ so that the exact solution to \cref{love} is $y(t)=1$. Since $r = 5$, this problem cannot be solved using Hale's method \cite{hal3} which works only for $r=1$. The resulting linear system constructed using the proposed approximation is an arrow-shaped one of infinite dimensions. We solve this system adaptively following the optimal truncation strategy of \cite{olv}. The computed solution is accurate in all but the last digit, as shown in \cref{fig:love}.

The proposed approximation allows the establishment of a spectral method for Fredholm convolution integral equations in the framework of infinite-dimensional linear algebra \cite{olv}.

%
%

\subsection{Pseudospectrum}\label{sec:convps}
As for solving Fredholm convolution integral equations, the calculation of the pseudospectra of the Fredholm convolution operator relies on the approximation to the operator. Consider the Huygens--Fresnel operator
\begin{align*}
(\mF_c^L u)(s) = \sqrt{\frac{iF}{\pi}}\int_{-1}^1 e^{-iF(s-t)^2}u(t)\md t,
\end{align*}
which comes from the modeling of the laser problem \cite[\S 60]{tre}. This integral operator acting on $u \in L_2([-1,1])$ correspond to stable resonators.

\begin{figure}[t!]
\centering
\includegraphics[width=0.5\linewidth]{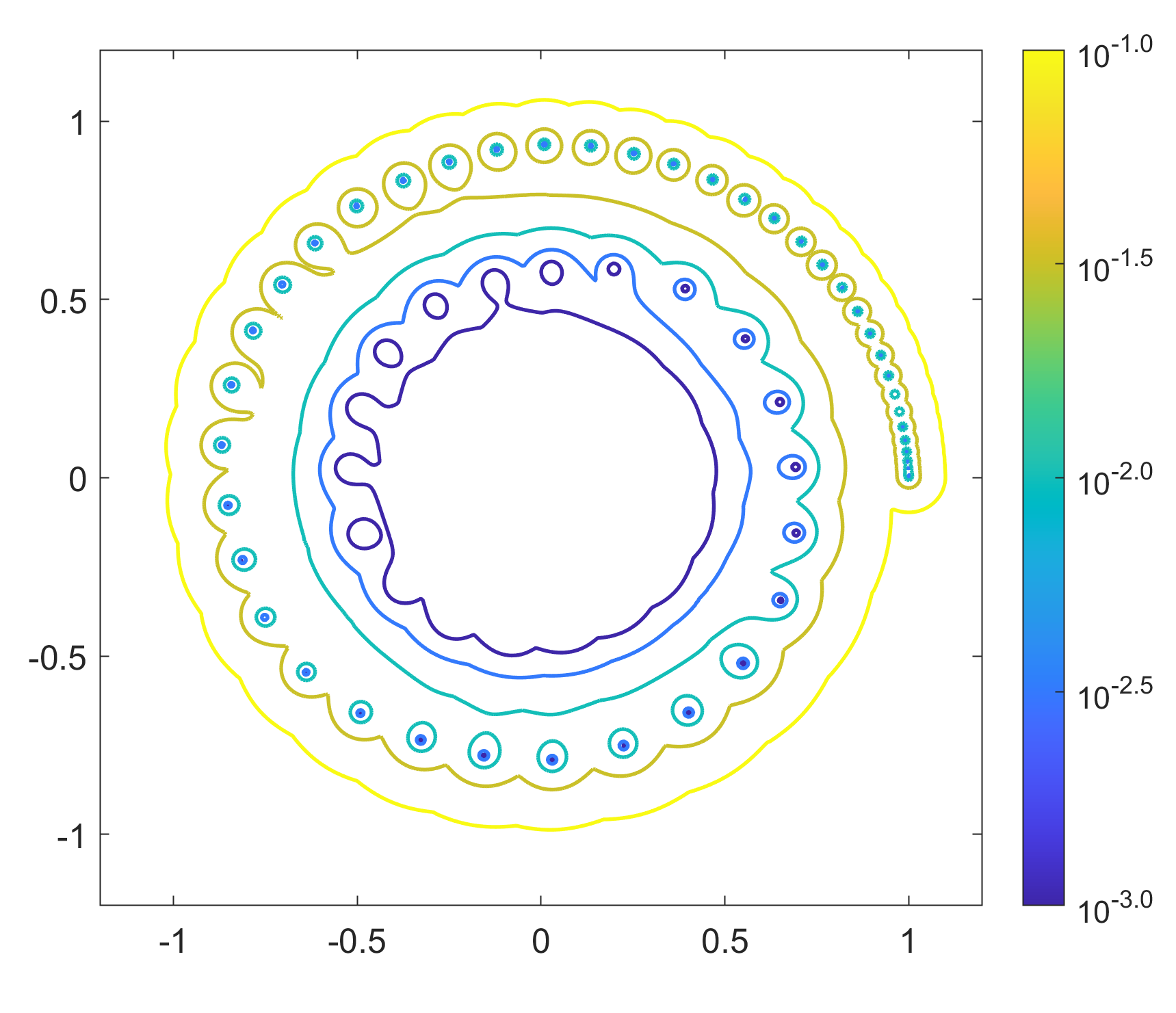}
\caption{$\varepsilon$-pseudospectra of the Huygens--Fresnel operator}
\label{fig:pseLaser}
\end{figure}

\cref{fig:pseLaser} displays the pseudospectra of the Huygens--Fresnel operator with the Fresnel number $F = 64\pi$. These boundaries of $\varepsilon$-pseudospectra at $\varepsilon = 10^{-1}, 10^{-1.5}, \dots, 10^{-3}$ are computed using the operator analogue of the Lanczos iteration \cite{den}, employing the ``solve-then-discretize'' strategy. For this operator-based method, the proposed spectral approximation plays an indispensable role. While the contours in \cref{fig:pseLaser} are visually indistinguishable from those in the final panel of Figure 60.2 from \cite[\S 60]{tre}, which Trefethen and Embree obtained following the conventional "discretize-then-solve" method, the traditional approach is at high risk of \textit{spectral pollution} and \textit{spectral invisibility} \cite{den}. It is only with the proposed spectral approximation we can, for the first time, obtain trustworthy pseudospectra of Fredholm convolution operators.

\section{Conclusion}\label{sec:conclusion}
We have presented a novel $\mO(M^2)$ method for constructing the spectral approximation of a convolution operator of Fredholm type defined by a Legendre series. In contrast to existing methods, the method is applicable to problems where the two convolved Legendre series are defined in intervals with an arbitrary length ratio. When such an approximation is used to evaluate the Fredholm convolution of two given Legendre series, the computational cost depends only on the length of the kernel, and is independent of the Legendre series being convolved and the interval length ratio. In all the cases that we test, the method is faster than the existing methods. In many cases where the Legendre series being convolved is of large degree or the interval length ratio is far greater than $1$, the new method is significantly faster. Moreover, such an approximation are of crucial importance in solving Fredholm convolution integral equations and the study of the spectral properties of the Fredholm convolution operators.

This work complements the previous work on the spectral approximation of the Volterra convolution operators \cite{hal2,xu2}. An exciting generalization of this work is to kernels with singularities following a similar approach in \cite{xu2} and this paper via recurrence relations. This may lead to new methods for the evaluation of fractional integrals and the solution to the fractional integral and differential equations. 



\bibliographystyle{siamplain}
\bibliography{fred}
\end{document}